\documentclass[11pt]{article}

\usepackage{amsmath, amsthm, amssymb}
\usepackage{graphicx, graphics}

\usepackage{color}
\def\bc{\begin{center}}
\def\ec{\end{center}}

\def\f{\frac}

\def\part#1#2{\f{\partial #1}{\partial #2}}

\newtheorem{Th}{Theorem}
\newtheorem{Lem}{Lemma}
\newtheorem{Prop}{Proposition}
\newtheorem{Def}{Definition}
\newtheorem{Cor}{Corollary}
\newtheorem{Rem}{Remark}

\def\pro{\smallskip\par\noindent{P\ r\ o\ o\ f: }}

\newcommand{\be}[1]{\begin{equation}\label{#1}}
\newcommand{\ee}{\end{equation}}
\newcommand{\beq} {\begin{eqnarray}}
\newcommand{\eeq} {\end{eqnarray}}

\newfont{\sss}{msam10 scaled\magstep0}
\newfont{\bbb}{msbm10 scaled\magstep1}
\newfont{\qqq}{msbm7}

\def\R{\mbox{\bbb R}}

\def\sq{{\sss \char3}}

\def\TTT{\end{document}}

\def\eps{\varepsilon}

\def\ga{\alpha}
\def\gb{\beta}

\def\gl{\lambda}
\def\gg{\gamma}

\def\gt{\tau}

\begin{document}
\bc
{ \Large \sc Stability by Lyapunov functions of Caputo fractional differential equations  with non-instantaneous impulses}
\ec

\vspace {1 cm} \bc {\large Ravi Agarwal$^{a,c,*}$, \hspace{0.2 cm} S. Hristova$^{b}$},  \hspace{0.2 cm}
D. O'Regan$^{d,c}$ \ec \vspace
{1 cm}

\begin{quote}
 $ ^{a}$ \ Department of Mathematics,\  Texas A$\&$M University-Kingsville,

 \hspace{0.3 cm} Kingsville,  TX 78363, \ USA

$^{b}$  Department of Applied Mathematics,\ Plovdiv University,  \  Plovdiv, \
Bulgaria

 $^{d}$ School of Mathematics, Statistics and Applied Mathematics,

 \hspace  {0.3 cm} National University of Ireland, Galway, Ireland

$^{c}$ NAAM Research Group, King Abdulaziz University,

 \hspace  {0.3 cm} Jeddah, Saudi Arabia

$^a$ \  {\it   e-mail  address}:  agarwal@tamuk.edu \ \ $^*$ corresponding author

 $^b$ \  {\it   e-mail  address}: snehri@gmail.bg

$^d$ \  {\it   e-mail  address}: donal.oregan@nuigalway.ie

\end{quote}

\vspace {1 cm}

\bc
{\sc Abstract}
\ec

$\ \ \ \ \  $ The stability of the zero solution of a nonlinear   Caputo fractional differential equation with noninstantaneous impulses is studied using
Lyapunov like functions.  The novelty of this paper is based on the  new  definition of the derivative of a Lyapunov like function along
the given noninstantaneous  impulsive fractional differential equations.  On one side this definition is a natural generalization of Caputo
 fractional Dini derivative of a function and on the other side it allows us the assumption for
Lyapunov functions to be weakened to continuity. By appropriate examples it
is shown the natural relationship between the defined derivative of Lyapunov
functions and Caputo derivative. 
  Several sufficient conditions for   uniform
stability and asymptotic uniform stability of the zero solution,  based on the new
definition of the derivative of Lyapunov functions  are established. Some examples are given
to illustrate the results. 

\vspace{0.5 cm}


    \vspace{0.5 cm}

\section{INTRODUCTION}\label{s1}

In the real world life  there are many processes and phenomena that are characterized by rapid changes in their state, so called impulses. Mainly there are two types of impulses:
\begin{itemize}
\item[-]  {\it instantaneous impulses}- the duration of these changes is relatively short compared to the overall duration of the whole process and the changes turn out to be irrelevant to the delopment of the studied process. The model of such kind of processes is given by impulsive differential equations (see, for example, \cite {AH}, \cite {H1}, \cite {H2}, \cite {H3}, \cite {H4}, \cite {H5}, \cite {BH1}, \cite {BH2}, \cite {HS},  \cite {E}, the monographs  \cite{S}, \cite {LB} and cited therein references);
\item[-] {\it noninstantaneous impulses} -  an impulsive action, which starts at an
arbitrary fixed point and keeps active on a finite time interval.  To achieve this aim E. Hernandez and  D. O'Regan  (\cite {HR}) introduced a new class of abstract differential equations for which the impulses are
not instantaneous and investigated the existence of mild and classical solutions. Later several authors studied and obtained qualitative properties of the solutions such as  existence (\cite {P}, \cite {PRR}, \cite {WL}, \cite {WZ}),  stability properties (\cite {SS}),  boundary value problems (\cite {LW}),  periodic solutions (\cite  {FW}). 
\end{itemize}

The stability of fractional order systems is quite recent. There are
several approaches in the literature to study stability, one of
which is the Lyapunov approach. Results on stability of fractional differential equations 
in the literature via Lyapunov functions could be divided into two
main groups:
\begin{itemize}
\item[-] continuously differentiable Lyapunov functions with their Caputo fractional  derivatives (see, for example, the
papers  \cite{DM1}, \cite{H},  \cite {Li1}, \cite{Li3}). 
\item[-] continuous Lyapunov functions with their fractional Dini derivatives (see, for example, the papers
\cite{devi vlm}, \cite {Lak fde}, \cite{LL}). 
\end{itemize}

Introducing impulses into fractional differential equations leads to
complications on the stability properties of the solutions. In the
literature there are several papers which examine existence and
qualitative properties of fractional differential equations with
instantaneous impulses (see, for example,  \cite {AB}, \cite{BK},
\cite {BS}, \cite {B},  \cite {ZN}). There are also results
concerning fractional differential equations with noninstantaneous
impulses; see for example \cite {ABT}, \cite {GD}, \cite {RD}, \cite
{RD1}, \cite {KP}, \cite {LX},  \cite {LW2}, \cite {11}, \cite {Y}.
In the above cited  papers the impulses start abruptly at some points  and their
action continues on a finite  interval. As a motivation for the
study of such kind of systems  we consider the following simplified
situation concerning the hemodynamical equilibrium of a person. In
the case of a decompensation (for example, high or low levels of
glucose) one can prescribe some intravenous drugs (insulin). Since
the introduction of the drugs in the bloodstream and the consequent
absorption for the body are gradual and continuous processes, we can
interpret the situation as an impulsive action which starts abruptly
and stays active on a finite time interval.

In this paper the stability of the zero solution of the
noninstantaneous impulsive  nonlinear
 Caputo fractional differential equations is studied. We
define in an appropriate way two types of fractional derivatives of Lyapunov functions: Fractional Dini derivative as well as Caputo fractional Dini derivative. Their applications is  discussed. Several sufficient conditions for  uniform stability and
asymptotic uniform stability are obtained. Some examples illustrating the obtained results are given.
 
\section{NOTES ON FRACTIONAL CALCULUS}\label{s2}
Fractional calculus generalizes the derivative and the integral of a
function to a non-integer order \cite{1, 2, Lak fde, podlubny, kilbas} and
there are several definitions of fractional derivatives and
fractional integrals.
In engineering, the fractional order $q$ is often less than 1, so
we restrict our attention to $q\in (0,1)$.

\textbf{1: } The Riemann--Liouville (RL)   fractional derivative   of order $q\in(0,1)$
of $m(t)$ is given by (see, for example, 1.4.1.1 \cite {1}, or \cite {podlubny})
\begin{equation*}
_{t_0}^{RL}D^{q}m(t)=\frac{1}{%
\Gamma \left( 1-q\right) }\frac{d}{dt}\int\limits_{t_0}^{t}\left(
t-s\right) ^{-q}m(s)ds,\ \ \ \  \   \  t\geq t_0.
\end{equation*}
where $\Gamma (.)$ denotes the Gamma
function.

\textbf{2: }  The Caputo fractional derivative of order $q\in(0,1)$
is defined by (see, for example, 1.4.1.3 \cite{1}) \be{198}
_{t_0}^{c}D^{q}m(t)=\frac{1}{\Gamma \left( 1-q\right)
}\int\limits_{t_0}^{t}\left( t-s\right) ^{-q}m^{\prime}(s)ds,\ \ \ \
\  t\geq t_0.  \ee The Caputo and Riemann-Liouville formulations
coincide when  $m(t_0)=0$. The properties of the Caputo derivative
are quite similar to those of ordinary derivatives. Also, the
initial conditions of fractional differential equations with the
Caputo derivative has a clear physical meaning and as a result the
Caputo derivative is usually used in real applications.

\textbf{3: } The Grunwald$-$Letnikov fractional derivative is given
by (see, for example, 1.4.1.2 \cite{1})
\begin{equation*}
_{t_0}^{GL}D^{q}m(t)=\lim_{h\to 0}\frac{1}{h^q}\sum_{r=0}^{[\frac{t-t_0}{h}]} (-1)^r (qCr)m(t-rh), \ \ \ \ t\geq t_0,
\end{equation*}and the Grunwald$-$Letnikov fractional Dini  derivative  by \be{765}
_{t_0}^{GL}D^{q}_{+}m(t)=\limsup _{h\to
0+}\frac{1}{h^q}\sum_{r=0}^{[\frac{t-t_0}{h}]} (-1)^r (qCr)m(t-rh),
\ \ \ \ t\geq t_0, \ee where $qCr$ are the Binomial coefficients and
$[\frac{t-t_0}{h}]$ denotes the integer part of the fraction
$\frac{t-t_0}{h}$.

The relations between the three types of fractional derivatives are
given by $_{t_0}^{c}D^{q}m(t)=\ _{t_0}^{RL}{D}^{q}[m(t)-m(t_0)]$

\begin{Prop}(Theorem 2.25 \cite {2}). Let $m\in C^1[t_0,b]$. Then, for $t \in (t_0,b]$, $_{t_0}^{GL}D^{q}m(t)=\
_{t_0}^{RL}D^{q}m(t)$ .
\end{Prop}

Also, according to Lemma 3.4 (\cite{2}), $\text{\
}_{t_0}^{c}D_{t}^{q}m(t)=_{t_0}^{RL}D_{t}^{q}m(t)-m(t_0)\frac{(t-t_0)^{-q}}{\Gamma(1-q)}$
holds.

 From the relation
between the Caputo fractional derivative and the Grunwald $-$ Letnikov
fractional derivative  using (\ref{765}) we
define the Caputo fractional Dini derivative  as
\begin{equation*}
_{t_0}^{c}D^{q}_+ m(t)=\ _{t_0}^{GL}D^{q}_{+}[m(t)-m(t_0)],
\end{equation*}
i.e.
\be{10}
_{t_0}^{c}D^{q}_+ m(t)=\limsup _{h\to 0+}\frac{1}{h^q}\Big [m(t)-m(t_0)-\sum_{r=1}^{[\frac{t-t_0}{h}]} (-1)^{r+1} (qCr)\Big ( m(t-rh)-m(t_0)\Big )\Big].
\ee

\begin{Def}(\cite{devi vlm}) We say $m\in C^q([t_0,T],\R^n)$ if  $m(t)$ is
differentiable (i.e. $m'(t)$ exists), the Caputo derivative $\ _{t_0}^{c}D^{q}m(t)$ exists  and
satisfies (\ref{198}) for $t\in[t_0,T]$.
\end{Def}

\begin{Rem}Definition 1 coud be extended to any interval $I\subset \R_+$. \end{Rem}

\begin {Rem} If  $m\in C^q([t_0,T],\R^n)$ then $_{t_0}^{c}D_+^{q}m(t)=\  _{t_0}^{c}D^{q}m(t)$.
\end{Rem}

\section{NONINSTANTANEOUS IMPULSES IN FRACTIONAL DIFFERENTIAL EQUATIONS}\label{s3}
In this paper we will assume two increasing sequences of points
$\{t_i\}_{i=1}^\infty$ and $\{s_i\}_{i=0}^\infty$ are given such
that $s_0=0<t_i\leq s_i<t_{i+1}$ , $i=1,2,\dots$, and $\lim_{k\to
\infty}t_k=\infty$. 

Let $t_0 \in \cup_{k=0}^{\infty}[s_k,t_{k+1})$ be a given arbitrary
point. Without loss of generality we will assume that $t_0\in
[s_0,t_{1})$, i.e. $ 0\leq t_0<t_1$. Define the sequence of points
$\{\gt_i\}_{i=0}^\infty$ by
\begin{displaymath}
   \gt_k= \left\{
     \begin{array}{lr}
        & t_0 \ \ \  \mbox {for}\ \  k=0, \\
&s_k\ \ \  \mbox {for}\ \  k\geq 1. \\
     \end{array}
   \right.
\end{displaymath}

Consider the initial value problem (IVP) for the system of \emph {noninstantaneous impulsive fractional differential equations} (NIFrDE)
with a Caputo derivative for $0<q<1$,%
\be{1}\begin{split}
& \text{\ }_{
\gt_k}^{c}D^{q}x=f(t,x) \ \mbox {for} \ t\in (\gt_k,t_{k+1}],\ k=0,1,2,\dots, \\
& x(t)=\phi_i(t,  x(t_i-0))\ \ \ \mbox {for} \  t\in (t_i,s_i], \ k=1,2,\dots, \\
& x(t_0)=x_0 \end{split} \ee  where $x,x_0\in \R^n$,  \ $f:\cup_{k=0}^{\infty} [\gt_k,t_{k+1}]\times \R^{n}\to \R^{n}$,\
  $\phi_i: [t_i,s_i]\times \R^n\to\R^n$,
$(i=1,2,3,\dots)$.

We suppose that the function $f(t,x)$ is smooth enough, such that for any initial  value $\tilde{x}_0\in
\R^n$ the IVP for the system of Caputo \emph{fractional differential equations} (FrDE)
\be{100} \text{\ }_{\gt_k}^{c}D^{q}x=f(t,x)  \ \ \  \mbox {for}\ \
t\in[\gt_k,t_{k+1}] \ \ \  \mbox {with}\ \  x(\gt_k)=\tilde{x}_0\ee
has a solution $x(t)=x(t;\gt_k,\tilde{x}_0)\in C^q([\gt_k,t_k],\R^n)$.
Some sufficient
conditions for global existence of solutions of (\ref{100}) are
given in \cite {BM}, \cite {Lak fde}.

 The IVP for FrDE (\ref{100}) is equivalent to the integral equation
\be{102}x(t)=\tilde{x}_0+\frac{1}{\Gamma(q)}\int_{\gt_k}^t\ (t-s)^{q-1}f(s,x(s))ds \ \ \  \mbox {for}\ \  t \in [\gt_k,t_{k+1}].\ee

\begin{Rem} The intervals $(t_k,s_k]$, $k=1,2,\dots$ are called intervals of noninstantaneous impulses and the functions $\phi_k(t,x)$, $k=1,2,\dots$, are called noninstantaneous impulsive functions.
\end{Rem}

 Now $t_0$ is the initial time. We will assume throughout this paper that the initial time  $t_0$ is not in an interval of noninstantaneous impulses, i.e. we will assume $t_0\in \cup_{k=0}^{\infty}[s_k,t_{k+1})$.

We will give a brief description of the solution of  IVP for NIFrDE (\ref{1}).
The solution $x(t;t_0,x_0)$, $t\geq t_0$ of (\ref{1}) is given by
\begin{displaymath}
   x(t;t_0,x_0)= \left\{
     \begin{array}{lr}
        & X_k(t) \ \ \  \mbox {for}\ \  t \in (s_{k},t_{k+1}],\ k=0,1,2,\dots, \\
&\phi_k(t,X_{k-1}(t_k-0))\ \ \  \mbox {for}\ \  t \in (t_k,s_k],\ \ k=1,2,\dots, \\
     \end{array}
   \right.
\end{displaymath}
where
\begin{itemize}
\item[-]$X_0(t)$ is the solution of IVP for FrDE (\ref{100}) for $k=0$, $t\in [t_0,t_{1}]$, $\tilde{x}_0=x_0$     and  $X_0(t)$ satisfies (\ref{102}) on  $ [t_0,t_{1}]$;
\item[-]$ X_{1}(t)$ is the solution of IVP for FrDE   (\ref{100}) for $k=1$, $t\in [s_{1},t_{2}]$, $\tilde{x}_0=\phi_{1}(s_{1}, X_0(t_{1}-0))$, and $X_{1}(t)$ satisfies (\ref{102}) on  $ [s_{1},t_{2}]$;
\item[-]$X_{2}(t)$ is the solution of IVP for  FrDE  (\ref{100}) for $k=2$, $t\in [s_{2},t_{3}]$, $\tilde{x}_0=\phi_{2}(s_{2}, X_{1}(t_{2}-0))$,  and  $X_{2}(t)$ satisfies (\ref{102}) on  $[s_{2},t_{3}]$;
\end{itemize}
and so on.

Also, the solution $x(t)=x(t;t_0,x_0)$, $t\geq t_0$ of (\ref{1}) is given by
\begin{displaymath}
   x(t)= \left\{
     \begin{array}{lr}
        & x_0+\frac{1}{\Gamma(q)}\int_{t_0}^t\ (t-s)^{q-1}f(s,x(s))ds \ \ \  \mbox {for}\ \  t \in [t_0,t_{1}],  \\
&\phi_1(t,x(t_1-0))\ \ \  \mbox {for}\ \  t \in (t_1,s_1],\\
& \phi_1(s_1,x(t_1-0))+\frac{1}{\Gamma(q)}\int_{s_1}^t\ (t-s)^{q-1}f(s,x(s))ds \ \ \  \mbox {for}\ \  t \in [s_1,t_{2}],  \\
&\phi_2(t,x(t_2-0))\ \ \  \mbox {for}\ \  t \in (t_2,s_2],\\
& \phi_2(s_2,x(t_2-0))+\frac{1}{\Gamma(q)}\int_{s_2}^t\ (t-s)^{q-1}f(s,x(s))ds \ \ \  \mbox {for}\ \  t \in [s_2,t_{3}],  \\
&..................................................................................\\
&\phi_k(t,x(t_k-0))\ \ \  \mbox {for}\ \  t \in (t_k,s_k],\\
& \phi_k(s_k,x(t_k-0))+\frac{1}{\Gamma(q)}\int_{s_k}^t\ (t-s)^{q-1}f(s,x(s))ds \ \ \  \mbox {for}\ \  t \in [s_k,t_{k+1}],  \\
&..................................................................................
     \end{array}
   \right.
\end{displaymath}

\begin{Rem} If $t_k=s_k$, $k=1,2,\dots $  then the IVP for NIFrDE (\ref{1}) reduces to an IVP for impulsive fractional differential equations  studied in  \cite {AB},   \cite{BK}, \cite {BS}, \cite {B}. In this case at any point of instantaneous impulse  $t_k$ the amount of jump of the solution $x(t)$ is given by  $\Delta x(t_k)=x(t_k+0)-x(t_k-0)=\Phi_k(x(t_k-0))=\phi_k(t_k,x(t_k-0))-x(t_k-0)$.
\end{Rem}

\begin{Rem} In the case $q=1$ the IVP for NIFrDE (\ref{1}) reduces to an IVP for noninstantaneous impulsive differential equations   studied in   \cite{LX},  \cite {SS},  \cite {LW}, \cite {WL}.
\end{Rem}

 \begin{Rem} In the case $q=1$,  $t_k=s_k$, $k=1,2,\dots $  the IVP for NIFrDE (\ref{1}) reduces to an IVP for
 impulsive  differential equations (see for example the books \cite{S}, \cite {LB} and the cited
references therein).
\end{Rem}

Let $J\subset \R_+$ be a given interval. We introduce the following
classes of functions \beq
PC^q(J)&=&\{u\in C^q(J\cap\big (\cup_{k=0}^{\infty}[s_k,t_{k+1})\big ),\R^n) \bigcup C(J\cap\big (\cup_{k=1}^{\infty}(t_k,s_{k}]\big ),\R^n):\nonumber \\
&& u(t_k)= u(t_{k}-0)=\lim_{t\uparrow t_{k}}u(t)<\infty,\ u(t_{k}+0)=\lim_{t\downarrow t_{k}}u(t)<\infty\nonumber\\
&&\ \ \ \ \ \ \ \ \ \ \ \ \ \ \ \ \ \ \ \ \ \ \ \ \ \ \ \ \ \ \ \mbox{for} \ k:\ t_k\in J,\nonumber \\
&& u(s_k)= u(s_{k}-0)=\lim_{t\uparrow s_{k}}u(t)=u(s_{k}+0)=\lim_{t\downarrow s_{k}}u(t)\ \ \mbox{for} \ k:\ s_k\in J\},\nonumber \\
PC(J)&=&\{u\in C(J\cap\big (\cup_{k=0}^{\infty}(t_k,t_{k+1})\big ),\R^n) :\nonumber \\
&& u(t_k)= u(t_{k}-0)=\lim_{t\uparrow t_{k}}u(t)<\infty\ \mbox{and}\ u(t_{k}+0)=\lim_{t\downarrow t_{k}}u(t)<\infty\nonumber \\
&&\ \ \ \ \ \ \ \ \ \ \ \ \ \ \ \ \ \ \ \ \ \ \ \ \ \ \ \ \ \ \  \mbox{for} \ k:\ t_k\in J\}.\nonumber \eeq

\begin{Rem} According to the above description any  solution of (\ref{1}) is from the
class $PC^q([t_0,b))$, $b\leq\infty$, i.e. any solution could have a
discontinuity at points $t_k,k=1,2,\dots$.
\end{Rem}

\textbf{Example 1}.   Consider  the IVP for the scalar NIFrDE
\be{311}\begin{split}
& \text{\ }_{\gt_k}^{c}D^{q}x=Ax \ \mbox {for} \ t\in (\gt_k,t_{k+1}],\ k=0,1,2,\dots, \\
&  x(t)=\Psi_k(t,x(t_k-0)) \ \mbox{for}\ \ \ t\in(t_k,s_k], \ k=1,2,\dots,\\
& x(t_0)=x_0, \end{split}
\ee
where $x, x_0\in\R$, $A$ is a constant.

The solution of (\ref{311}) is given by
\begin{displaymath}
   x(t;0,x_0)= \left\{
     \begin{array}{lr}
 &x_0E_q(A(t-t_0)^q)  \ \ \ \ \ \ \ \mbox {for}\ \   t \in[t_0,t_1]  \ \ \ \ \ \ \ \ \ \ \ \ \ \ \ \ \ \ \ \ \ \ \ \ \ \\
&\Psi_k(t, x(t_k-0)) \ \ \ \ \ \  \ \ \ \    \mbox {for}\ \  t \in (t_k,s_k], \ \ k=1,2,\dots  \ \ \ \ \ \\
    & \Psi_k\Big(s_k, x(t_k-0)\Big)E_q\Big(A(t-s_k)^q\Big)  \ \ \ \ \ \ \ \mbox {for}\ \  t \in [s_k,t_{k+1}],\ k=1,2,3,\dots
     \end{array}
   \right.
\end{displaymath}  where the
Mittag $-$ Leffler function (with   one parameter) is defined by
$E_{q}(z)=\sum\limits_{k=0}^{\infty }\frac{z^{k}}{\Gamma (qk+1)}.$

If $\Psi_k(t,x)=a_k(t)x, \ a_k:\ [t_k,s_k]\to\R, \ k=1,2,3,\dots$,
then the solution of NIFrDE (\ref {311}) is given by
\be{107}
\begin{split}
 x(t;0,x_0)=
\begin{cases}
x_0E_q(A(t-t_0)^q)\ \ \ \ \ \ \ \ \ \ \ \ \ \ \ \ \ \ \text{for}\ t \in[0,t_1],  \\
 x_0E_q(AC^q) \Big(\prod _{i=1}^{k-1}a_i(s_i)E_q(AC_i^q)\Big) a_k(t)& \\
  \ \ \ \ \ \ \ \ \ \ \ \ \ \ \ \ \ \ \text{for}\  t \in (t_k,s_k], \ k= 1,2,\dots, &\\
x_0E_q(AC^q)a_k(s_k)\Big(\prod _{i=1}^{k-1}a_i(s_i)E_q(AC_i^q)\Big)E_q(A(t-s_k)^q) \\
\ \ \ \ \ \ \ \ \ \ \ \ \ \ \ \ \ \  \text{for}\ t \in (s_k,t_{k+1}], \ k=1,2,\dots. &\nonumber
     \end{cases}
\end{split}
\ee
where $C=t_1-t_0$ and $C_k=s_k-t_k\geq 0$, $k=1,2,\dots$.

If $A=0$ and $\Psi_k(t,x)=a_k(t)x, \ a_k:\ [t_k,s_k]\to\R, \
k=1,2,3,\dots$  the solution of NIFrDE (\ref {311}) is given by
\begin{displaymath}
   x(t;0,x_0)= \left\{
     \begin{array}{lr}
        &x_0  \ \ \ \ \ \ \ \mbox {for}\ \   t \in[t_0,t_1], \ \ \ \ \ \ \ \ \ \ \ \ \ \ \ \ \   \\
        &x_0 \Big(\prod _{i=1}^{k-1}a_i(s_i)\Big) a_k(t)\ \ \ \ \mbox {for}\ \  t \in (t_k,s_k], \ k=1,2,\dots\\
 &x_0\Big(\prod _{i=1}^{k}a_i(s_i)\Big)\ \ \ \ \mbox {for}\ \  t \in (s_k,t_{k+1}], \ k=1,2,\dots
     \end{array}
   \right.
\end{displaymath}

\hfill~~~\sq

\section{DEFINITIONS CONCERNING STABILITY AND LYAPUNOV FUNCTIONS}\label{s4}

The goal of the paper is to study the stability properties of the
system NIFrDEs (\ref{1}). In the definition below we denote by
$x(t;t_{0},x_{0})\in PC^q([t_0,\infty),\R^n)$ any solution of
(\ref{1}).
 \begin{Def}
  The zero solution
of the IVP for NIFrDE (\ref{1}) is said to be
\begin{itemize}
 \item \emph {stable } if for every $\epsilon >0$  and $t_0\in \cup_{k=0}^{\infty}[s_k,t_{k+1})$ there
 exist  $\delta =\delta (\epsilon,t_0)>0$ such that
 for any $x_0\in  \R^n$
   the inequality  $||x_0||<\delta$  implies $||x(t;t_0,x_0)||<\epsilon$ for $t\geq t_0$;
 \item \emph {uniformly stable } if for every $\epsilon >0$   there
 exist  $\delta =\delta (\epsilon)>0$ such that
 for any initial point $t_0\in \cup_{k=0}^{\infty}[s_k,t_{k+1})$ and any inital value $x_0\in \R^n$ with $||x_0||<\delta$  the inequality $||x(t;t_0,x_0)||<\epsilon$ holds for $t\geq t_0$;
 \item  \emph { uniformly attractive} if for $\gb>0:$  for
 every $\epsilon >0$   there exist $T=T(\epsilon)>0$ such that for any initial pont $t_0\in \cup_{k=0}^{\infty}[s_k,t_{k+1})$ and any initial value $x_0\in \R^n$ with
 $||x_0||<\gb$ the inequality $||x(t;t_0,x_0)||<\epsilon$ holds for $t\geq t_0+T$;
 \item  \emph {uniformly asymptotically stable} if the zero solution is uniformly stable  and uniformly attractive.
 \end{itemize}
  \end{Def}

\textbf{Example 2}. Consider the scalar NIFrDE  (\ref{311}) where
$A\leq  0$ and $\Psi_k(t,x)=a_k(t)x$ and  $a_k:\ [t_k,s_k]\to\R,  \
k=1,2,3,\dots$  are such that $\sup_{t\in[t_k,s_k]}|a_k(t)|\leq
M_k$, $\prod _{i=1}^{\infty}M_i<\infty$ where $M_k>0$ are constants.
According to Example 1 and the inequality $0<E_{q}(A(T-\gt)^{q})\leq
1$ for $T\geq \gt$
  there exists a constant $M>0$ such that
\be{421} \left\vert x(t;t_0,x_0) \right\vert \leq M \left\vert
x_{0}\right\vert\ \ \ \ \mbox{for}\ \ t\geq t_0. \ee Inequality (\ref{421}) guarantees that the zero
solution of (\ref{311}) is uniformly stable.

\hfill~~~\sq

In this paper we will use the followings sets:
\begin{eqnarray*}
\mathcal{K} &\mathcal{=}&\{a\in C[\mathbb{R}_{+},\mathbb{\!R}_{+}]:a\text{
is strictly increasing and }a(0)=0\}, \\
S(A ) &=&\{x\in \R^n:||x||\leq A \}, \ \ \ \ A>0.
\end{eqnarray*}

We now introduce the class $\Lambda $ of Lyapunov-like functions
which will be used to investigate the  stability of the zero solution of the
system NIFrDE (\ref{1}).

\begin{Def} Let $J\in  \mathbb{\ R}_+$ be a given interval, and $\Delta \subset
\mathbb{\ R}^n,\ \ 0\in \Delta$ be a given set. We will say that the function $V(t,x):J\times
\Delta\rightarrow \mathbb{R}_{+}, $ $%
V(t,0)\equiv 0$ belongs to the class $\Lambda(J,\Delta) $ if
\begin{itemize}
\item[1.]The function $%
V(t,x)$ is continuous on $J/\{t_k\in J\}\times \Delta$  and it  is  locally Lipschitzian with respect to its second
argument;
\item[2.]For each $t_k\in Int(J)$ and $x\in \Delta$ there exist finite limits
$$V(t_k-0,x)=\lim_{t\uparrow t_k}V(t,x)<\infty,\ \ \ \ \mbox{and}\ \  V(t_k+0,x)=\lim_{t\downarrow t_k}V(t,x)<\infty$$
and the following equalities are valid
$$V(t_k-0,x)=V(t_k,x).$$
\end{itemize}
\end{Def}
\begin{Rem} In the case when the Lyapunov function does not depend on the time $t$, i.e. $V(t,x)=V(x)\in C[\Delta,\R_+], $ $%
V(0)= 0$,
 $\Delta \subset
\mathbb{\ R}^n,\ \ 0\in \Delta$, and the function $%
V(x)$ is  locally Lipschitzian then we will say  $V(x)\in \Lambda^C(\Delta).$
\end{Rem}
Lyapunov-like functions used to discuss stability for differential
equations require an appropriate definition of their derivatives along the studied differential equations.  For
fractional differential equations some authors (see, for example,
 \cite{Li1}, \cite{Li3}) used the so called Caputo fractional
derivative of Lyapunov function $\text{ }_{t_0}^{c}D^{q}V(t,x(t))$ where $x(t)$ is the unknown solution of the studied fractional differential equation.  This approach requires the
function to be smooth enough (at least continuously differentiable)
and also some conditions involved are quite restrictive. Other
authors used the so called Dini fractional  derivative of Lyapunov
function (\cite {Lak fde}, \cite{LL}). This is based on the Dini derivative of the Lyapunov function $V(t,x)$ among the ordinary differential equation  $x'=f(t,x)$ given by
\be{3331} DV(t,x)=\limsup_{h\to
0}\frac{1}{h}\Big[V(t,x)-V(t-h,x-hf(t,x)\Big]. \ee
The authors generalized (\ref{3331}) to the {\emph {
Dini fractional  derivative}} along the FrDE $ \ _{t_0}^{c}{ D}^{q}x=f(t,x),\ t\geq t_0$  by
\be{444} ^{c}{ D}_{+}^{q}V(t,x)=\limsup_{h\to
0}\frac{1}{h^q}\Big[V(t,x)-V(t-h,x-h^qf(t,x)\Big]. \ee
This  definition requires
only the continuity of the Lyapunov function.

In this paper we  will use piecewise continuous Lyapunov functions
from the above introduced class $\Lambda ([t_0,T),\Delta)$. We will
introduce the derivative of  Lyapunov function in two different ways
and we will discuss their applications.

We now  define the {\emph {generalized Caputo fractional Dini
derivative}} of the Lyapunov-like function $V(t,x)\in \Lambda ([t_0,T),\Delta)$
along trajectories of solutions of IVP for the system NIFrDE
(\ref{1}). It is based on the Caputo fractional Dini derivative of a
function $m(t)$ given by (\ref {10}).  The Caputo fractional Dini
derivative along trajectories of solutions of IVP for the system
NIFrDE (\ref{1}) is given by: \be{20}
\begin{split}
&\text{ }_{(\ref{1})}^{c}D_{+}^{q}V(t,x;t_0,x_0) =\limsup_{h\rightarrow 0^{+}} {\frac{1}{h^{q}}}\bigg\{
V(t,x)-V(t_0,x_0)\\
&\ \ \ \  -\sum_{r=1}^{[\frac{t-t_0}{h}]}(-1)^{r+1}qCr \bigg [ V(t-rh,x-h^{q}f(t,x))-V(t_0,x_0)\bigg ]\bigg\}\\
&\ \ \ \ \ \ \ \ \ \ \ \  \ \ \mbox{for}\  t\in (s_k,t_{k+1})\cap (t_0,T),\ k=0,1,2,\dots,
\end{split}
\ee where  $x,x_0\in \Delta$, and for any $t\in (s_k,t_{k+1})\cap
(t_0,T)$ there exists $h_t>0$ such that $t-h\in (s_k,t_{k+1})\cap
(t_0,T)$, $x-h^{q}f(t,x)\in \Delta $ for $ 0<h\leq h_t$.

\begin{Rem} The generalized Caputo fractional Dini  derivative of the Lyapunov function  was
introduced and used for studying stability properties of the zero
solution of Caputo fractional differential equations in \cite {AHR}.
\end{Rem}

The formula (\ref{20}) could be reduced to \be{200}
\begin{split}
&\text{ }_{(\ref{1})}^{c}D_{+}^{q}V(t,x;t_0,x_0)  \\
&\ \ \ \ \ \ \ =\limsup_{h\rightarrow 0^{+}} {\frac{1}{h^{q}}}\bigg\{
V(t,x)-\sum_{r=1}^{[\frac{t-t_0}{h}]}(-1)^{r+1}qCr V(t-rh,x-h^{q}f(t,x))\bigg\}\\
&\ \ \ \ \ \ \ \ \ \ \ \ \ \ -V(t_0,x_0) \frac{(t-t_0)^{-q}}{\Gamma(1-q)}  \ \ \mbox{for}\  t\in (s_k,t_{k+1})\cap (t_0,T),\ k=0,1,2,\dots,
\end{split}
\ee

Now,  based on the  Dini fractional derivative of a continuous
Lyapunov function defined by (\ref{444}), we will define  the {\emph {generalized Dini fractional derivative}}
of the Lyapunov-like function $V(t,x)\in \Lambda ([t_0,T),\Delta)$ along
trajectories of solutions of IVP for the system NIFrDE (\ref{1})
 by
\be{14}\begin{split}
&\text{ }_{(\ref{1})}^{c}{\mathcal D}_{+}^{q}V(t,x)=\lim_{h\rightarrow 0^{+}}\sup {\frac{1}{h^{q}}}%
\Big[V(t,x)-V(t-h,x-h^{q}f(t,x))\Big]\\
&\ \ \ \ \ \ \ \ \ \ \ \ \ \  \ \ \ \  \ \ \ \ \ \ \ \  \mbox{for}\  t\in (s_k,t_{k+1})\cap (t_0,T),\ k=0,1,2,\dots ,
\end{split}\ee
where  $x\in \Delta$, and for any $t\in (s_k,t_{k+1})\cap (t_0,T)$
there exists $h_t>0$ such that $t-h\in (s_k,t_{k+1})\cap (t_0,T)$,
$x-h^{q}f(t,x)\in \Delta $ for $ 0<h\leq h_t$.

\textbf{Example 3}.  Let $V\in \Lambda(\R_+,\R)$ be given by 
$V(t,x)=m(t)g(x)$  where the function $m\in
C^1(\cup_{i=0}^\infty(s_k,t_{k+1}),\R_+)$, $g:\R^n\to \R_+$ is a locally Lipshitz function  such that the limit  $$\ ^{Fr}D_qg(x)=\lim_{h\rightarrow 0^{+}}\sup \frac{g(x)-g(x-h^{q}f(t,x))}{h^q}$$ exists for $x\in \R^n$.

 First we apply the formula
(\ref{14}) to obtain the generalized  Dini fractional derivative of the considered Lyapunov function. We obtain
\be{148}\begin{split}
&\text{ }_{(\ref{1})}^{c}{\mathcal D}_{+}^{q}V(t,x) =\lim_{h\rightarrow 0^{+}}\sup {\frac{1}{h^{q}}}\Big[m(t)g(x)-m(t-h)g\big(x-h^{q}f(t,x)\big)\Big]
 \\
 &=m(t)\lim_{h\rightarrow 0^{+}}\sup \frac{g(x)-g(x-h^{q}f(t,x))}{h^q}\\
&\ \ \ \ \ \ +\Big(\lim_{h\rightarrow 0^{+}}\sup \frac{m(t)-m(t-h)}{h}\Big)\Big(\lim_{h\rightarrow 0^{+}}\sup h^{1-q}g\big(x-h^{q}f(t,x)\big)\Big)\\
 &= m(t)\ ^{Fr}D_qg(x)\ \ \ \ \mbox{for}\ \ \ t\in (s_k,t_{k+1}),\ k=0,1,2,\dots.
\end{split}\ee

Next we use  (\ref{200})  to obtain the generalized Caputo fractional Dini derivative of the function $V$.  Let $t\in (s_k,t_{k+1}),\ k=0,1,2,\dots$. Apply  equalities (\ref{765}),  $\
_{t_0}^{RL}D^{q}1=\frac{(t-t_0)^{-q}}{\Gamma(1-q)}$,  and 
$\lim_{h\rightarrow 0^{+}}\sup {\frac{1}{h^{q}}}\sum\limits_{r=0}^{[\frac{t-t_0}{%
h}]}\left( -1\right) ^{r}qCr\ m(t-rh)=\
_{t_0}^{RL}D^{q}\Big(m(t)\Big)$  to (\ref{200}) and  obtain
 \be{1777}
\begin{split}
&\text{ }_{(\ref{1})}^{c}D_{+}^{q}V(t,x;t_0, x_0) \\
&=\lim_{h\rightarrow 0^{+}}\sup {\frac{1}{h^{q}}} \Big[ m(t)g(x) +g( x-h^{q}f(t,x))\sum\limits_{r=1}^{[\frac{t-t_0}{h}]}\left( -1\right) ^{r}qCr\ m(t-rh)\Big]\\
&\ \ \ \ \ \ -m(t_0)g(x_0)\frac{(t-t_0)^{-q}}{\Gamma(1-q)}\\
&=m(t)\lim_{h\rightarrow 0^{+}}\sup {\frac{g(x)-g( x-h^{q}f(t,x))}{h^{q}}} \\
&\ \ \ +\lim_{h\rightarrow 0}g( x-h^{q}f(t,x))\lim_{h\rightarrow 0^{+}}\sup {\frac{1}{h^{q}}}\sum\limits_{r=0}^{[\frac{t-t_0}{h}]}\left( -1\right) ^{r}qCr\ m(t-rh)\Big]\\
&\ \ \ \ \ \ -m(t_0)g(x_0)\frac{(t-t_0)^{-q}}{\Gamma(1-q)}\\
&=m(t)\ ^{Fr}D_qg(x) +g(x)\ _{t_0}^{RL}D^{q}\Big(m(t)\Big)
- m(t_0)g(x_0)\frac{(t-t_0)^{-q}}{\Gamma(1-q)}, \\
&\ \ \ \ \ \ \ \ \ \mbox{for}\ \ \ \ \ t\in (s_k,t_{k+1}),\ k=0,1,2,\dots.
\end{split}\ee
or
\be{173}
\begin{split}
&\text{ }_{(\ref{1})}^{c}D_{+}^{q}V(t,x;t_0, x_0)  \\
&=m(t)\ ^{Fr}D_qg(x) +g(x)\ _{t_0}^{C}D^{q}\Big(m(t)\Big)
+\Big (g(x)-g(x_{0})\Big)m(t_0)\frac{(t-t_0)^{-q}}{\Gamma(1-q)}, \\
&\ \ \ \ \ \ \ \ \ \mbox{for}\ \ \ \ \ t\in (s_k,t_{k+1}),\ k=0,1,2,\dots.
\end{split}\ee

\textbf{Example 4}.  Let $V\in \Lambda(\R_+,\R)$ be given by  the equality 
$V(t,x)=m(t)x^2$  where $m\in
C^1(\cup_{i=0}^\infty(s_k,t_{k+1}),\R_+)$  and  $x\in \R$.

In this  case \be {777}\begin{split}
\ ^{Fr}D_qg(x)&=\ ^{Fr}D_q(x^2)=\lim_{h\rightarrow 0^{+}}\sup \frac{x^2-(x-h^{q}f(t,x))^2}{h^q}\\
&=\lim_{h\rightarrow 0} f(t,x)(2x-h^{q}f(t,x))=2xf(t,x)\end{split}\ee

First we apply the formula (\ref{14}) to obtain the generalized Dini
fractional derivative of the considered Lyapunov function. From (\ref{148}) and (\ref{777}) we obtain 
\be{16}\begin{split}
&\text{ }_{(\ref{1})}^{c}{\mathcal D}_{+}^{q}V(t,x) =2 x\ m(t)f(t,x),\ \ \ \ \mbox{for}\ \ \ t\in (s_k,t_{k+1}),\ k=0,1,2,\dots.
\end{split}\ee
Note $\text{ }_{(\ref{1})}^{c}{\mathcal D}_{+}^{q}V(t,x)$ does not
depend on the order $q$ of the fractional differential equation. The
behavior of solutions of fractional differential equations depends
significantly on the order $q$. For example, let us consider the
simple fractional differential equation $ \text{\
}_{0}^{c}D^{q}x+x(t)=1,\ x(0)=0$ whose solution is given by
$x(t)=t^q E_{q,1+q}(-t^q)$. From Figure 1 it can be seen $lim_{t\to
\infty}x(t) =a$ where $a$ is different for different values of the
order $q$ of fractional differential equation.

\begin{figure}
\begin{center}
\begin{minipage}[b]{0.45\linewidth}
\centering
\includegraphics[width=1\textwidth]{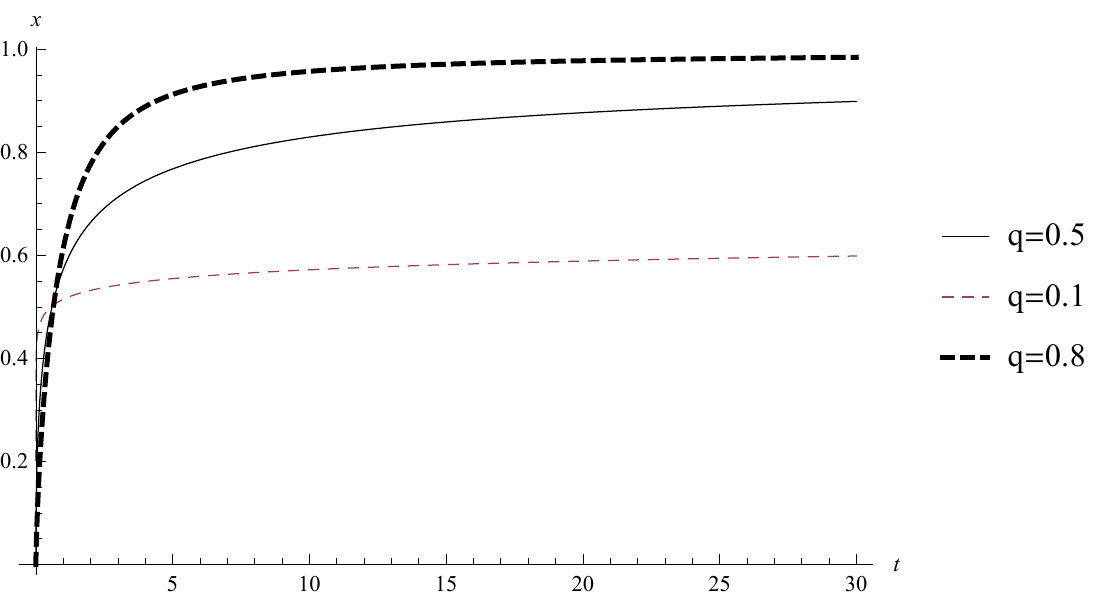}
{\small \it Figure 1.  Graphs of solutions  of $\text{\ }_{0}^{c}D^{q}x+x(t)=1$ for various $q$. }
\end{minipage}\hfill
\begin{minipage}[b]{0.45\linewidth}
\centering
\includegraphics[width=1\textwidth]{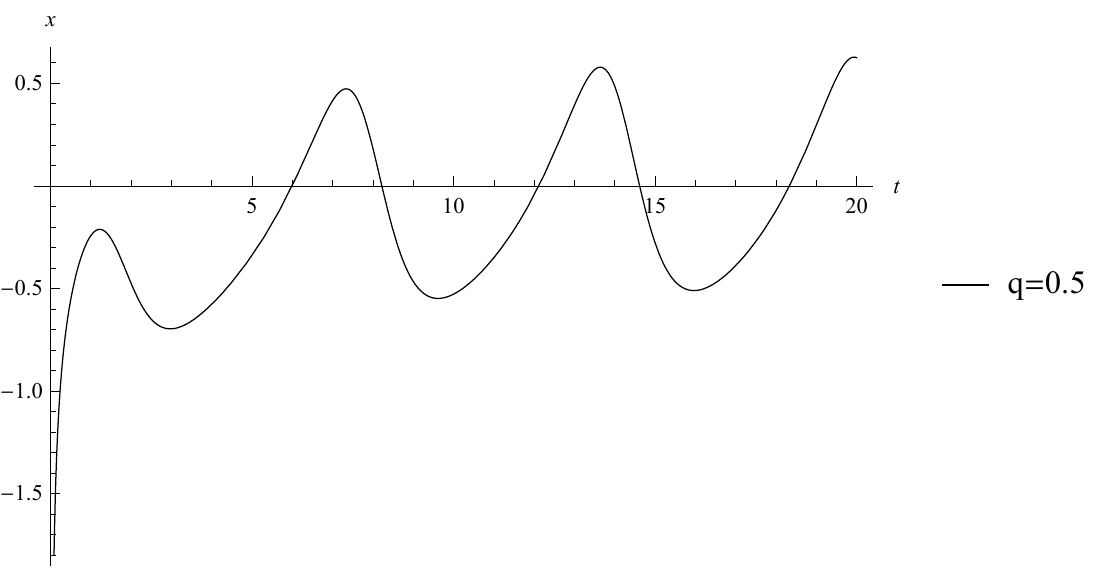}
{\small \it Figure 2.  Example 8. Graph of  $f(t)$. }
\end{minipage}
\end{center}
 \end{figure}

 Next we use  (\ref{200})  to obtain the generalized Caputo fractional Dini derivative of the function $V$.  Let $t\in (s_k,t_{k+1}),\ k=0,1,2,\dots$. From (\ref{1777}), (\ref{173}) for $g(x)=x^2$ and (\ref{777}) we obtain
 \be{177}
\begin{split}
&\text{ }_{(\ref{1})}^{c}D_{+}^{q}V(t,x;t_0, x_0) =\\
&=2x\ m^2(t)f(t,x) +x^2\ _{t_0}^{RL}D^{q}\Big(m^2(t)\Big)
-(x_{0})^2 m^2(t_0)\frac{(t-t_0)^{-q}}{\Gamma(1-q)}, \\
&\ \ \ \ \ \ \ \ \ \mbox{for}\ \ \ \ \ t\in (s_k,t_{k+1}),\ k=0,1,2,\dots.
\end{split}\ee
or
\be{17}
\begin{split}
&\text{ }_{(\ref{1})}^{c}D_{+}^{q}V(t,x;t_0, x_0)  \\
&=2x\ m^2(t)f(t,x) +x^2\ _{t_0}^{C}D^{q}\Big(m^2(t)\Big)
+\Big (x^2-(x_{0})^2 \Big)m^2(t_0)\frac{(t-t_0)^{-q}}{\Gamma(1-q)}, \\
&\ \ \ \ \ \ \ \ \ \mbox{for}\ \ \ \ \ t\in (s_k,t_{k+1}),\ k=0,1,2,\dots.
\end{split}\ee
Note the generalized Caputo fractional Dini derivative $\text{
}_{(\ref{1})}^{c}D_{+}^{q}V(t,x;t_0, x_0) $ depends significantly
not only on the order $q$ of the fractional differential equation
but also on the initial data.

 The derivative of the Lyapunov function in the
well known case for first order impulsive  differential equations
($q=1$) is \be{18} D_{+}V(t,x)=2x\ m^2(t)f(t,x)+x^{2}
\frac{d}{dt}\Big[m^2(t)\Big],\ \ \ t\in (s_k,t_{k+1}),\
k=0,1,2,\dots. \ee

The generalized Caputo fractional Dini derivative given by formula (\ref{200})
seems to be the  natural generalization of the derivative of
Lyapunov functions  for ordinary  differential equations (with or
without impulses).

 \hfill~~~\sq

\section{COMPARISON RESULTS}\label{s5}
Again in this section we assume $0<q<1$. Also, $t_0\in
\cup_{k=0}^{\infty}[s_k,t_{k+1})$ so without loss of generality we
assume $t_0\in [s_0,t_{1})$. We will obtain some comparison results
for NIFrDE (\ref{1}) using the definitions (\ref{200}) and
(\ref{14}) for a derivative of Lyapunov-like function.

\subsection{Generalized Caputo fractional Dini derivative}\label{s5.1}
 In this section we will use the following result for  fractional differential equations:

 \begin {Lem}  (\cite {AHR}) Assume the following conditions are satisfied:
 \begin{itemize}
\item[1.] The function $x^*(t)=x(t;t_0,x_0)\in C^q([t_0,T],\Delta),$ is a solution of the FrDE
\be{1001} \text{ }_{t_0}^{c}D^{q}x=f(t,x),\ \ \ t\in [ t_0,T]\ \ \
\mbox{with}\ x(t_0)=x_0 \ee  where $\Delta \subset \R^n,\
0\in\Delta$,  $x_0\in\Delta$ and  $t_0,\ T\in\R_+,\, t_0<T$ are
given constants.
\item[2.]  The function $V\in \Lambda^C([t_0,T],\Delta)$ and
   the inequality $$\text{ }_{(\ref{1001})}^{c}D_{+}^{q}V(t,x;t_0,x_0 )\leq 0\ \ \mbox{for}\ \ (t,x)\in [t_0,T]\times \Delta$$
   holds.
 \end{itemize}

   Then the inequality    $V(t,x^*(t))\leq V(t_0,x_0)$ holds  for $t\in [t_0,T].$
   \end {Lem}

Now we  will prove some  comparison results for noninstantaneous
impulsive Caputo fractional differential equations.

\begin {Lem} ({\it Comparison result for NIFrDE by generalized Caputo fractional Dini derivative}). Let:
\begin{itemize}
\item[  1.] The function $x^*(t)=x(t;t_0,x_0)\in PC^q([t_0,T],\Delta)$ is a solution of
the NIFrDE (\ref{1}) where $\Delta \subset \R^n,\ 0\in\Delta$, $x_0\in\Delta$ and
 $t_0,T$ are given constants such that $t_0\in [s_0,t_{1})$, $T>t_0$.
  \item[2.]  The function $V\in \Lambda([t_0,T],\Delta)$ and
    \begin{itemize}
\item[(i)] the inequality 
\be{1678}\text{ }_{(\ref{1})}^{c}D_{+}^{q}V(t,x^*(t);t_0,x_0 )\leq 0\ \ \ \mbox{for}\ \ t\in (t_0,T)\bigcap \cup_{k=0}^{\infty}(s_k, t_{k+1})\ee
  holds;
\item[(ii)]     the inequalities
$$V(t,x^*(t))\leq V(t_k-0,x^*(t_k-0))\ \ \mbox{for}\ t\in [t_0,T]\bigcap(t_k,s_k]\ \mbox{for}\ \ k=1,2,3,\dots$$
 hold.
\end{itemize}
\end{itemize}
   Then the inequality     $V(t,x^*(t))\leq V(t_0,x_0)$ holds on  $[t_0,T]$.
   \end {Lem}

    \pro   We  use induction to prove Lemma 2.

Let $t\in[t_0,t_{1}]\cap[t_0,T]$. The function $x^*(t)\in
C^q([t_0,t_{1}]\cap[t_0,T],\R^n)$, satisfies  FrDE (\ref{1001}) and
from Lemma  1  (with  $T=t_{1}$) the inequality $V(t,x^*(t))\leq
V(t_0,x_0)$ holds on  $[t_0,t_{1}]\cap[t_0,T]$.

Let $T>t_{1}$ and $t\in(t_{1},s_{1}]\cap[t_0,T]$.  From condition
2(ii) and the above we get  $V(t,x^*(t))\leq
V(t_{1}-0,x^*(t_{1}-0))=V(t_{1},x^*(t_{1})) \leq V(t_0,x_0)$.

Let $T>s_{1}$ and $t\in (s_{1},t_{2}]\cap[t_0,T]$. Consider  the
function $\overline{x}_{1}(t)= x^*(t)$ for $t\in (s_{1},t_{2}]$  and
$\overline{x}_{1}(s_{1})= x^*(s_{1})=\phi_{1}(s_{1}, x^*(t_{1}-0))$.
The function  $\overline{x}_{1}(t)\in C^q([s_{1},t_{2}],\R^n)$ and
satisfies IVP for FrDE (\ref{1001}) with $t_0=s_{1}$,  $x_0=
x^*(s_{1})$,  and $T=t_{2}$.   Using condition 2(i),  Lemma 1   for
the function $\overline{x}_{1}(t)$ and the above we obtain
$V(t,x^*(t))=V(t,\overline{x}_{1}(t))\leq
V(s_{1},\overline{x}_{1}(s_{1}))= V(s_{1},x^*(s_{1}))\leq
V(t_0,x_0)$.

Continue this process and an induction argument proves the claim of Lemma 2 is
true for $t\in [t_0,T]$.

\hfill~~~\sq

     \begin {Lem} (\cite {AHR}).  Let the following conditions be satisfied:
        \begin{itemize}
        \item[1.] The function $x^*(t)=x(t;t_0,x_0)\in C^q([t_0,T],\Delta)$  is  a solution
        of the FrDE (\ref{1001})
        where $\Delta \subset \R^n,\ 0\in\Delta$,  $x_0\in\Delta$, $t_0,\ T\in\R_+,\ t_0<T$ are given constants.
        \item[2.]
        The function $V\in \Lambda^C([t_0,T],\Delta)$ is such that  for any points $t \in [t_0,T]$, $x\in \Delta$
   the inequality $$\text{ }_{(\ref{100})}^{c}D_{+}^{q}V(t,x ;\gt_0,x_0)\leq -c(||x||)$$
   holds where $c\in \mathcal{K}$.
\end{itemize}

   Then for $t\in [t_0,T]$ the inequality  \be{456}V(t,x^*(t))\leq V(t_0,x_0)-\frac{1}{\Gamma(q)}\int_{t_0}^t (t-s)^{q-1}c(||x^*(s)||)ds   \ee holds.
   \end {Lem}

\begin {Lem} ({\it Comparison result for NIFrDE, negative generalized Caputo fractional Dini derivtive}).  Assume the following conditions are satisfied:
\begin{itemize}
\item[  1.] The function $x^*(t)=x(t;t_0,x_0)\in PC^q([t_0,T],\Delta)$ is a solution of
the NIFrDE (\ref{1}) where $\Delta \subset \R^n,\ 0\in\Delta$, $x_0\in\Delta$ and
 $t_0,T$ are given constants such that $t_0\in [s_0,t_1)$, $T>t_0$.
  \item[2.]  The function $V\in \Lambda([t_0,T],\Delta)$ and
    \begin{itemize}
\item[(i)] the inequality $\text{ }_{(\ref{1})}^{c}D_{+}^{q}V(t,x^*(t);t_0,x_0 )\leq -c(||x^*(t)||)$   for $t\in (t_0,T)\bigcap \cup_{k=0}^{\infty}(s_k, t_{k+1})$
  holds where $c\in \mathcal{K}$;
\item[(ii)] for any  $k=1,2\dots$    the inequalities
$$V(t,x^*(t))\leq V(t_k-0,x^*(t_k-0))\ \ \mbox{for}\ t\in [t_0,T]\bigcap(t_k,s_k]$$  hold.
\end{itemize}
\end{itemize}

   Then for  $t\in [t_0,T]$ the inequality
    \begin{displaymath}
   V(t,x^*(t))\leq  \left\{
     \begin{array}{lr}
    &V(t_0,x_0)-\frac{1}{\Gamma(q)}\int_{t_0}^t (t-s)^{q-1}c(||x^*(s)||)ds,\ \ \ \ \ \ \ \ \ \ \ \ \ \ \ \ \ \ \ \ \ \ \ \ \ \ \ \ \ \ \ \ \ \  \\
		& \  t \in [t_0,t_{1}] \ \ \ \ \ \ \ \ \ \ \ \ \ \ \ \ \ \ \ \ \ \ \ \ \ \ \ \ \ \  \\
 & V(t_0,x_0)-\frac{1}{\Gamma(q)}\Big(\sum_{i=0}^{k-1}\int_{\gt_i}^{t_{i+1}} (t_{i+1}-s)^{q-1}c(||x^*(s)||)ds\  \ \ \ \ \ \ \ \ \ \ \ \ \ \ \ \ \ \ \ \ \ \ \ \ \ \ \  \ \ \ \\
&\ \ \ \ \  +\int_{s_k}^{t} (t-s)^{q-1}c(||x^*(s)||)ds\Big), \ \ \ \ \ \ \ \ \ \ \ \ \ \ \ \ \ \ \ \ \ \ \ \ \ \ \ \ \ \ \ \ \ \  \\
&  t \in (s_k,t_{k+1}]\cap(t_0,T), \ k\geq 1 \ \ \ \ \ \ \ \ \ \ \ \ \ \ \ \ \ \\
&V(t_0,x_0)-\frac{1}{\Gamma(q)}\sum_{i=0}^{k-1}\int_{\gt_i}^{t_{i+1}} (t_{i+1}-s)^{q-1}c(||x^*(s)||)ds,\ \ \ \ \ \ \ \ \ \ \ \ \ \ \ \ \ \ \ \ \ \ \ \ \ \ \ \ \ \ \ \ \ \ \\
& \  t \in (t_{k},s_{k}]\cap[t_0,T],\ k\geq 1\ \ \ \ \ \ \ \ \ \ \ \ \ \ \ \ \ 
     \end{array}
   \right.
\end{displaymath}
    holds  where
    \begin{displaymath}
   \gt_k= \left\{
     \begin{array}{lr}
    &t_0\ \ \mbox{for}\ \ k=0, \\
 & s_k\ \ \mbox{for}\ \ k\geq 1.
     \end{array}
   \right.
\end{displaymath}
   \end {Lem}
    \pro Let $t\in[t_0,t_{1}]\cap[t_0,T]$. The function $x^*(t)\in C^q([t_0,t_1]\cap[t_0,T],\Delta)$ and
    satisfies the IVP for FrDE (\ref{1001})  for  $T=t_1$.  From Lemma 3   the inequality $V(t,x^*(t))\leq V(t_0,x_0)-\frac{1}{\Gamma(q)}\int_{t_0}^t (t-s)^{q-1}c(||x^*(s)||)ds$ holds, i.e. the claim of Lemma 4 is true  on $[t_0,t_{1}]\cap[t_0,T]$.

Let $T>t_{1}$ and $t\in(t_{1},s_{1}]\cap[t_0,T]$.  From condition
2(ii) and the above we get
\be{567}\begin{split}V(t,x^*(t))&\leq  V(t_{1}-0,x^*(t_{1}-0))=V(t_{1},x^*(t_{1})) \\
& \leq V(t_0,x_0)-\frac{1}{\Gamma(q)}\int_{t_0}^{t_{1}} (t_{1}-s)^{q-1}c(||x^*(s)||)ds.\nonumber \end{split}\ee

Let  $T>s_{1}$ and   $t\in(s_{1},t_{2}]\cap[t_0,T]$. Consider  the
function $\overline{x}_{1}(t)= x^*(t)$ for $t\in (s_{1},t_{2}]$  and
$\overline{x}_{1}(s_{1})= x^*(s_{1})=\phi_{1}(s_{1}, x^*(t_{1}-0))$.
The function  $\overline{x}_{1}(t)\in C^q([s_{1},t_{2}],\R^n)$ and
satisfies IVP for FrDE (\ref{1001})   with $t_0=s_{1}$,  $x_0=
x^*(s_{1})$ and $T=t_{2}$.   Using condition 2(i),  Lemma 3   for
the  function $\overline{x}_{1}(t)$, and the above we obtain \beq
&&V(t,x^*(t))=V(t,\overline{x}_{1}(t))\leq V(s_{1}+0,\overline{x}_{1}(s_{1}))-\frac{1}{\Gamma(q)}\int_{s_{1}}^t (t-s)^{q-1}c(||x^*(s)||)ds\nonumber\\
&& = V(s_{1},x^*(s_{1}))-\frac{1}{\Gamma(q)}\int_{s_{1}}^t (t-s)^{q-1}c(||x^*(s)||)ds
\nonumber\\
&& \leq V(t_0,x_0)-\frac{1}{\Gamma(q)}\int_{t_0}^{t_{1}} (t_{1}-s)^{q-1}c(||x^*(s)||)ds-\frac{1}{\Gamma(q)}\int_{s_{1}}^t (t-s)^{q-1}c(||x^*(s)||)ds.
\nonumber
\eeq
Therefore, the claim of Lemma 4 is true on $(s_{1},t_{2}]\cap[t_0,T]$.

Let  $T>t_{2}$ and    $t\in(t_{2},s_{2}]\cap[t_0,T]$.  From
condition 2(ii) and the above we obtain
\be{5671}\begin{split} &V(t,x^*(t))\leq V(t_{2}-0,x^*(t_{2}-0))= V(t_{2},x^*(t_{2}))\\
&\leq V(t_0,x_0)-\frac{1}{\Gamma(q)}\int_{t_0}^{t_{1}} (t_{1}-s)^{q-1}c(||x^*(s)||)ds -\frac{1}{\Gamma(q)}\int_{s_{1}}^{t_{2}} (t_{2}-s)^{q-1}c(||x^*(s)||)ds.\nonumber \end{split}\ee

Continue this process and an induction argument proves the claim is
true for $t\in [t_0,T]$.

\hfill~~~\sq

\begin{Rem}
The results of Lemma 2 and Lemma 4 are true on the half line (recall
\cite {AHR} that Lemma 1 and Lemma 3 extend to the half line).
\end{Rem}

\begin{Rem} The results of Lemma 2 and  Lemma 4  will be similar with slight changes of condition 2(ii) if the initial time $t_0$ is in a interval of noninstantaneous  impulses, i.e. $t_0\in \cup_{k=1}^{\infty}(t_k,s_k]$.
\end{Rem}

\subsection{Generalized Dini fractional derivative}
Now we present the analogue results of Section 5.1 using
(\ref{444}).

\begin {Lem} (\cite {LL}). Let $V\in C(\R_+\times \R^n , \R_+)$ and $V(t,x)$ be locally Lipschitzian in $x$. Assume that  $$^cD_+^qV(t,x)\leq g(t,V(t,x)),\ \ \  (t,x)\in \R_+\times \R^n,$$ where $g\in C(\R_+^2,\R)$ and $^cD_+^qV(t,x)$ is defined by (\ref{444}). Suppose that the maximal solution $r(t;t_0,u_0)$ of IVP  $$ \text{\ }_{t_0}^{c}D^{q}u=g(t,u),\ \ \ u(t_0)=u_0\geq 0$$ exists on $[t_0,\infty)$. Then $V(t_0,u_0)\leq u_0$ implies $V(t,x(t))\leq r(t)$ for $t\geq t_0$ where $x(t)=x(t;t_0,x_0)$ is any solution of IVP (\ref{1001})  existing on $[t_0,\infty)$.
\end{Lem}

\begin{Cor}Let $V\in C(\R_+\times\R^n , \R_+)$ and $V(t,x)$ be locally Lipschitz in $x$ and $\ ^cD_+^qV(t,x)\leq 0$ for $(t,x)\in \R_+\times \R^n.$   Then  $V(t,x(t))\leq V(t_0,x_0)$ for $t\geq t_0$ where $x(t)=x(t;t_0,x_0)$ is any solution of IVP (\ref{1001})  existing on $[t_0,\infty)$.
\end{Cor}

Now we will give a comparison result for NIFrDE (\ref{1}) by application of generalized Dini fractional derivative of Lyapunov function. 
\begin {Lem} ({\it Comparison result for NIFrDE by generalized Dini fractional derivative}).  Assume the conditions of Lemma 2 are satisfied where the inequality (\ref{1678}) is replaced by 
\be{689}\text{ }_{(\ref{1})} ^c{\mathcal D}_+^qV(t,x^*(t))\leq0\ \ \ \mbox{for}\ \ t\in [t_0,\infty)\bigcap\cup_{k=0}^{\infty}(s_k, t_{k+1}).\ee
 
   Then the inequality     $V(t,x^*(t))\leq V(t_0,x_0)$ holds on  $[t_0,\infty)$.
   \end {Lem}
        \pro The proof of Lemma 6 is similar to that in Lemma 2 where instead of Lemma 1 we apply Corollary 1.

\hfill~~~\sq

If $V(x)=x^Tx=\sum_{k=1}^n x_k^2\in \Lambda^C(\R^n),\ \
x=(x_1,x_2,\dots,x_n)^T$ we obtain the following result:
\begin {Cor} ({\it Comparison result by a quadratic Lyapunov function}). Let the function
$x^*(t)=x(t;t_0,x_0)\in PC^q([t_0,\infty),\R^n)$ be a solution of the NIFrDE (\ref{1})  such that
    \begin{itemize}
\item[(i)]  $\big(x^*(t)\big)^T f(t,x^*(t))\leq 0$  for $t\in [t_0,\infty)\bigcap (s_k, t_{k+1}),\ \ k=0,1,2,\dots$.
\item[(ii)] for any  $k=1,2\dots$    the inequalities
$$||x^*(t)||\leq ||x^*(\gt_k-0)||\ \ \mbox{for}\ t\in [t_0,\infty)\bigcap(t_k,s_k]\ \mbox{with}\ \ \gt_k=max\{t_0,t_k\}$$
  hold.
\end{itemize}
   Then the inequality     $||x^*(t)||\leq ||x_0||$ holds on  $[t_0,\infty)$.
   \end {Cor}
    \pro The proof follows from  Lemma 6 applied with $V(x)=x^Tx=\sum_{k=1}^n x_k^2$  and
    note
    \be{161}\begin{split}
&\text{ }_{(\ref{1})}^{c}{\mathcal D}_{+}^{q}V(x)
=\lim_{h\rightarrow 0^{+}}\sup {\frac{1}{h^{q}}}\Big[\sum_{i=1}^n
x_i^2-\sum_{i=1}^n \Big(x_i-h^{q}f^{(i)}(t,x)\Big)^2\Big]
 \\
 &=\sum_{i=1}^n2x_i f^{(i)}(t,x)=2\ x^Tf(t,x)\\
&\ \ \ \ \ \ \ \ \ \ \mbox{for}\ \ \ t\in [t_0,\infty)\bigcap (s_k,t_{k+1}),\ k=0,1,2,\dots
\end{split}\ee

     \hfill~~~\sq

\textbf{Example 5}.   Consider the scalar NIFrDE  (\ref{311}) where
$A<0$, and $\Psi_k(t,x)=a_k(t)x$, and $a_k:\ [t_{k},s_k]\to\R,  \
k=1,2,3,\dots$  are such that $\sup_{t\in[t_k,s_k]}|a_k(t)|\leq 1$.
Its solution $x^*(t)=x(t;0,x_0)$ is given in Example 1.  Using
$0<E_{q}(A(T-\gt)^{q})\leq 1$ for $T\geq \gt$ it follows that the
conditions of Corollary 2 are satisfied  and therefore $|x^*(t)|\leq
|x_0|$ for $t\geq 0$.

    \hfill~~~\sq

\begin{Rem} In this paper we assumed an infinite number of points
$\{t_i\}_{i=1}^{\infty}$, $\{s_i\}_{i=1}^{\infty}$ with $0<t_i \leq
s_i<t_{i+1}$ and $\lim_{k\to \infty}t_k=\infty$. However it is worth
noting that the results in Section 5 (and elsewhere) hold true if we
only consider a finite of points $\{t_i\}_{i=1}^{p}$,
$\{s_i\}_{i=1}^{p}$ with $0<t_i \leq s_i<t_{i+1}$, $i=1,....,p$ with
$t_{p+1}=T$.
\end{Rem}

\section{MAIN RESULTS}\label{s6}
We will obtain sufficient conditions for
stability of the zero solution of nonlinear impulsive Caputo fractional differential equations.  Again
we assume $0<q<1$.

We  say conditions (H) are satisfied if :

{\bf (H1)} The function $f\in C(\cup_{k=0}^\infty[s_k,t_{k+1}],\R^n)$,  $f(t,0)\equiv 0$ for $t\in \cup_{k=0}^\infty[s_k,t_{k+1}]$  is  such that  for any initial point $(\tilde{t}_0,\tilde{x}_0)\in\cup_{k=0}^\infty[s_k,t_{k+1})\times\R^n$, $s_k\leq \tilde{t}_0<t_{p+1}$, $p$ is a nonzero integer,  the IVP for the system of FrDE (\ref{100}) with $\gt_k=\tilde{t}_0$ has a solution $x(t;\tilde{t}_0,\tilde{x}_0)\in C^q([\tilde{t}_0,t_{p+1}],\R^n)$.

{\bf (H2)} The functions $\phi_k \in C([t_k,s_k]\times \R^n,\R^n)$, $\phi_k(t, 0)\equiv  0$ for $t\in[t_k,s_k]$, $k=1,2,\dots$.

\begin{Rem} Conditions (H) guarantee the  existence of a
solution $x(t;t_0,x_0)\in PC^q([t_0,\infty),R^n)$ of NIFrDE (\ref{1})
for any initial data $(t_0,x_0)\in \cup_{k=0}^\infty[s_k,t_{k+1})\times\R^n.$ If (H) is
satisfied then (\ref{1}) has a zero solution. \end{Rem}

\subsection{Stability by generalized Caputo fractional Dini derivative}\label{s6.1}
\begin{Th} ({\it Stability}). Let the following conditions be satisfied:
\begin{itemize}
\item[1.] Condition (H) is satisfied.
\item[2.] There exists a function $V\in \Lambda (\R_+, \R^n) $  such that $V(t,0)= 0$  and
\begin{itemize}
\item[(i)] the inequality
\be{13517}\text{ }_{(\ref{1})}^{c}D_{+}^{q}V(t,x ;t_0,x_0)\leq 0  \ \ \ \ \ \ \ \mbox{for}\ \  t\in \cup_{k=0}^{\infty}(s_k,t_{k+1}),\ \  x, x_0\in \R^n \ee holds;
\item[(ii)]  for any point $x\in \R^n$ and any $t\in (t_k,s_k]$, $k=1,2,3,\dots$ the inequality $V(t,\phi_k(t,x))\leq V(t_k-0,x)$  holds;
\item[(iii)] $b(||x||)\leq V(t,x)$ for $ t\in\R_+, \ x\in\R^n ,$
  where  $b \in \mathcal{K}$.
\end{itemize}
\end{itemize}
     Then  the zero solution of the NIFrDE (\ref{1}) is   stable.
     \end{Th}

\pro  Let $\epsilon>0$ and $t_0\in \cup_{k=0}^{\infty}[s_k,t_{k+1})$
be arbitrary  given numbers. There exists a $p\in \{0,1,...\}$ with
$t_0\in [s_p,t_{p+1})$.  Without loss of generality assume $p=0$.

Since $V(t_0,0)= 0$  there exists $\delta_1=\delta_1(t_0,\eps)>0$
such that $V(t_0,x)<b(\eps)$ for $||x||<\delta_1$.  Let  $x_0\in
\R^n$ with   $ ||x_0||<\delta_1$. Then $V(t_0,x_0)<b(\eps)$.
Consider any solution $x^*(t)=x(t;t_0,x_0)\in
PC^q([t_0,\infty),\R^n)$ of  NIFrDE (\ref {1}).  From inequality
(\ref{13517}) it follows that $$\text{
}_{(\ref{1})}^{c}D_{+}^{q}V(t,x^*(t);t_0,x_0)\leq 0 \ \ \mbox{for}\
t\in (t_0,\infty)\bigcap \cup_{k=0}^{\infty}(s_k,t_{k+1}),$$ i.e.
condition 2(i) of Lemma 2 (with $T=\infty$, see Remark 10) is
satisfied.

Let $t\in (t_k,s_k]\bigcap[t_0,\infty)$, $ k=1,2,3,\dots$. From
condition 2(ii) of Theorem 1 we get
$$V(t,x^*(t))=V(t,\phi_k(t,x^*(t_k)))\leq  V(t_k-0,x^*(t_k))=
V(t_k-0,x^*(t_k-0)).$$ Therefore, condition 2(ii) of Lemma 2 is
fulfilled.

From Lemma 2 applied to the solution $x^*(t)$ with $T=\infty$ (see
Remark 10) and condition 2($iii$)  we obtain
$$ b(||x^*(t)||)\leq V(t,x^*(t))\leq V(t_0,x_0)<b(\epsilon),$$
so the result follows.

 \hfill~~~\sq

  \begin{Th}  ({\it Uniform stability}). Let the following conditions be satisfied:
\begin{itemize}
\item[1.] Condition (H) is satisfied.
\item[2.] There exists a function $V\in \Lambda (\R_+,\R^n) $  such that
    \begin{itemize}
    \item[(i)] the inequality is satisfied
\be{135}\text{ }_{(\ref{1})}^{c}D_{+}^{q}V(t,x ;t_0,x_0)\leq 0  \ \ \ \ \ \ \ \mbox{for}\ \  t\in \cup_{k=0}^{\infty}(s_k,t_{k+1}),\ \  x, x_0\in S(\lambda) \ee where $\gl>0$ is a given number;
\item[(ii)]  for any point  $t\in (t_k,s_k]$, $k=1,2,3,\dots$ and any $x\in S(\gl)$  the inequality $V(t,\phi_k(t,x))\leq V(t_k-0,x)$  holds;
\item[(iii)]  $b(||x||)\leq V(t,x)\leq a(||x||)$ for $ t\in\R_+, \ x\in \R^n,$
  where  $a,b \in \mathcal{K}$.
\end{itemize}
\end{itemize}
     Then  the  zero solution of NIFrDE (\ref{1})  is  uniformly stable.
\end {Th}

\pro Let $\epsilon\in (0,\gl]$  and $t_0\in
\cup_{k=0}^{\infty}[s_k,t_{k+1})$  be arbitrary  given numbers.
There exists  $p\in \{0,1,...\}$ with $t_0\in [s_p,t_{p+1})$.
Without loss of generality assume $p=0$.

Let $\delta_1<
\min\{ \epsilon, b(\epsilon)\}$. From $a \in \mathcal{K}$ there
exists  $\delta_2 =\delta_2 (\epsilon)>0 $ so if $s<\delta_2$ then
$a(s)<\delta_1$.  Let $\delta=\min (\epsilon, \delta_2)$. Choose the
initial value  $x_0\in \R^n$ such  that $||x_0||<\delta$ and let
$x^*(t)=x(t;t_0,x_0)\in PC^q([t_0,\infty),\R^n)$ be a solution of the IVP for NIFrDE (\ref
{1}). We now prove that \be{9011}||x^*(t)||<\epsilon,\ \ \ \ t\geq
t_0.\ee

Assume inequality (\ref{9011}) is not true and let $t^*=inf\{t>t_0:\ ||x^*(t)||\geq \epsilon\}$. Then \be{555}||x^*(t)||<\epsilon \ \ \ \mbox{for}\ t\in[t_0,t^*)\ \ \ \mbox{and} \ ||x^*(t^*)||=\epsilon .\ee
If  $ t ^* \not = t_k,\ k=1,2,\dots$ or if $t^*=t_p$ for some natural number $p$ and $||x^*(t_p-0)||=\epsilon$ then (\ref{555}) is true. If for  a natural number $p$ we have $t^*=t_p$ and $||x^*(t_p-0)||<\epsilon$, then according to Lemma 2 for $ T=t_p$ and $\Delta=S(\gl)$ we obtain $V(t,x^*(t))\leq V(t_0,x_0)$ for $t\in[t_0,t_p]$. Then for all $t\in(t_p,s_p]$ from condition 2(iii) we get $b(||x^*(t)||)\leq V(t,x^*(t))=V(t,\phi_p(t,x^*(t_p-0)))\leq V(t_p-0,x^*(t_p-0)) \leq V(t_0,x_0)\leq
a(\delta)<\delta_1<b(\eps)$. Thus $||x^*(t)| | < b^{-1}(\delta_1)< \epsilon$ for $t \in (t_p,s_p]$,  and this contradicts the choice of $t^*$. Therefore, (\ref{555}) holds.

Then, $x^*(t)\in S(\gl)$ on $[t_0,t^*]$ and conditions 2(i) and
2(ii) of Lemma 2 are satisfied on $[t_0,t^*]$. From Lemma 2 applied
to the solution $ x^*(t)$ with $T=t^*$ and $\Delta=S(\gl)$ we get
$V(t,x^*(t))\leq V(t_0,x_0)$  on $[t_0,t^*]$. Then applying
condition 2 (iii) of Theorem 2 we obtain
$b(\eps)=b(||x^*(t^*)||)\leq V(t^*,x^*(t^*))\leq V(t_0,x_0)\leq
a(\delta)<\delta_1<b(\eps)$. The contradiction proves (\ref {9011})
and therefore, the zero solution of  NIFrDE (\ref{1}) is uniformly
stable.

\hfill~~~\sq

Now we present some sufficient conditions for the uniform asymptotic
stability of the zero solution of the NIFrDE.

\begin{Th} ({\it Uniform asymptotic stability}). Let the following conditions be satisfied:
\begin{itemize}
\item[1.] Condition (H) is satisfied.
\item[2.] There exists a positive constant $M<\infty$  such that $\sum_{i=1}^{\infty}(s_{i}-t_i)\leq M.$
 \item[3.] There exists a function $V\in \Lambda (\R_+,\R^n) $  such that
    \begin{itemize}
    \item[(i)] the inequality is satisfied
\be{1351}\text{ }_{(\ref{1})}^{c}D_{+}^{q}V(t,x ;t_0,x_0)\leq -c(||x||) \ \ \ \ \ \ \ \mbox{for}\ \  t\in \cup_{k=0}^{\infty}(s_k,t_{k+1}),\ \  x, x_0\in S(\lambda) \ee where $\gl>0$ is a given number, $c\in \mathcal{K}$;
\item[(ii)]  for any point  $t\in (t_k,s_k]$, $k=1,2,3,\dots$ and any $x\in S(\gl)$  the inequality $V(t,\phi_k(t,x))\leq V(t_k-0,x)$  holds;
\item[(iii)]  $b(||x||)\leq V(t,x)\leq a(||x||)$ for $ t\in\R_+, \ x\in \R^n,$
  where  $a,b \in \mathcal{K}$.
\end{itemize}
\end{itemize}

     Then  the zero solution  of NIFrDE (\ref{1}) is uniformly asymptotically stable.
\end {Th}

\pro  From Theorem 2 the zero solution of the NIFrDE (\ref{1}) is
uniformly stable. Therefore, for the number $\gl$  there exists
$\ga=\ga(\gl)\in (0,\gl)$ such that for any $\tilde{t}_0 \in \cup_{k=0}^{\infty}[s_k,t_{k+1})$
 and
$\tilde{x}_0  \in\R^n$   the inequality  $||\tilde{x}_0||<\ga$
implies  \be{789}||x(t;\tilde{t}_0,\tilde{x}_0)||<\gl\ \ \
\mbox{for}\ \ t\geq \tilde{t}_0\ee where
$x(t;\tilde{t}_0,\tilde{x}_0)$ is any solution of the NIFrDE (\ref{1})
(with initial data $(\tilde{t}_0,\tilde{x}_0)$).

Now we will prove that the zero solution of the fractional
differential equations (\ref{1}) is uniformly attractive. Consider
the constant $\gb\in(0,\ga]$  such that $a(\gb)\leq b(\ga)$. Let
$\epsilon \in(0,\gl]$  and $t_0\in \cup_{k=0}^{\infty}[s_k,t_{k+1})$
be arbitrary  given numbers. There exists a $p\in \{0,1,...\}$ with
$t_0\in [s_p,t_{p+1})$.  Without loss of generality assume $p=0$.

Let the point $x_0\in\R^n,\, ||x_0||<\gb$ and  $x^*(t)=x(t;t_0,x_0)$
be any solution of (\ref{1}). Then $b(||x_0||)\leq
a(||x_0||)<a(\gb)<b(\ga)$, i.e.  $||x_0||<\ga$ and according to
(\ref{789})  the   inequality \be{88} ||x^*(t)||<\gl\ \ \
\mbox{for}\ \ t\geq t_0 \ee holds, i.e. the solution $x^*(t)\in
S(\gl)$  on $[t_0,\infty)$.

Choose a constant $\gg=\gg(\epsilon)\in
(0,\epsilon]$ such that $a(\gg)<b(\epsilon)$. Let  $T>\sqrt[q]{a(\ga)\frac{q\Gamma(q)}{c(\gg)}}+M$
and  $m$ be a natural number such that $s_m<t_0+T\leq t_{m+1}$. Note $T$ depends only on $\eps$ but not on $t_0$. We now prove
that \be{4450} ||x^*(t)||< \epsilon\ \ \mbox {for }\ t\geq t_0+T.
\ee Assume \be{4490} ||x^*(t)||\geq \gg\ \ \mbox {for every}\ t\in
[t_0,t_0+T].\ee Then from Lemma 4 (applied to  the interval
$[t_0,t_0+T]$ and $\Delta=S(\gl)$), conditions 2 and  3 (ii) of Theorem 3, inequality $a^q+b^q\geq (a+b)^q$ for $a,b >0$  and the choice of $T$ we get
\be{4419}\begin{split}
& V(t_0+T,x^*(t_0+T))\\
& \leq V(t_0,x_0)-\frac{1}{\Gamma(q)}\Big(\sum_{i=0}^{m-1}\int_{\gt_i}^{t_{i+1}} (t_{i+1}-s)^{q-1}c(||x^*(s)||)ds\\
&\ \ \ \ \ \ \ \ \ \ \ \ \ \ \ \ +\int_{s_m}^{t_0+T} (t_0+T-s)^{q-1}c(||x^*(s)||)ds\Big)\\
& \leq a(||x_0||)-\frac{c(\gg)}{\Gamma(q)}\Big (\sum_{i=0}^{m-1}\int_{\gt_i}^{t_{i+1}} (t_{i+1}-s)^{q-1}ds+\int_{s_m}^{t_0+T} (t_0+T-s)^{q-1}ds\Big)\\
&<a(\ga)-\frac{c(\gg)}{q\Gamma(q)}\Big ((t_1-t_0)^q+\sum_{i=1}^{m-1}(t_{i+1}-s_i)^{q} +(T+t_0-s_m)^{q} \Big) \\
&\leq a(\ga)-\frac{c(\gg)}{q\Gamma(q)}\Big ((t_1-t_0)+\sum_{i=1}^{m-1}(t_{i+1}-s_i) +(T+t_0-s_m) \Big)^q   \\
&= a(\ga)-\frac{c(\gg)}{q\Gamma(q)}\Big (-\sum_{i=1}^{m}(s_{i}-t_i)+T \Big)^q \leq a(\ga)-\frac{c(\gg)}{q\Gamma(q)}\Big (-M+T \Big)^q  <0.\nonumber
\end{split}
\ee

The above contradiction proves there exists $t^*\in [t_0,t_0+T]$ such
that $||x^*(t^*)||<\gg$.  Let the natural number $p$ be such that $t_{p-1}\leq t^* < t_p$.

{\it Case 1.} Let $t\in [t^*,t_p]$.

If $s_{p-1}<t^*<t_p$ then for $t\in [t^*,t_p]$ the function $x^*(t)\in C^q([t^*,t_p],\R^n)$
 and from Lemma 3 we get $V(t,x^*(t))\leq V(t^*,x^*(t^*))-\frac{1}{\Gamma(q)}\int_{t^*}^{t} (t-s)^{q-1}c(||x^*(s)||)ds\leq V(t^*,x^*(t^*))$.

If $t_{p-1}<t^*\leq s_{p-1}$ then for $t\in [t^*,t_p]$ the function
$x^*(t)\in PC^q([t^*,t_p],\R^n)$ and from Lemma 4 we get
$V(t,x^*(t))\leq V(t^*,x^*(t^*))$.

{\it Case 2.}  For  any  $t> t^*, \ t \in (s_k,t_{k+1}],\ k= p, p+1,
\dots,  $ from Lemma 4 for $\Delta=S(\gl)$  we obtain
\be{3342}\begin{split}
&V(t,x^*(t))\\
&\leq V(t^*,x^*(t^*))-\frac{1}{\Gamma(q)}\Big(\int_{t^*}^{t_{p}} (t-s)^{q-1}c(||x^*(s)||)ds+\sum_{i=p}^{k-1}\int_{s_i}^{t_{i+1}} (t_{i+1}-s)^{q-1}c(||x^*(s)||)ds\\
&\ \ \ \ \ \ \ \ \ \ \ \ \ \ +\int_{s_k}^{t} (t-s)^{q-1}c(||x^*(s)||)ds\Big)\leq V(t^*,x^*(t^*)). \nonumber
\end{split}
\ee

{\it Case 3.}  For  any $t> t^*, \ t \in (t_{k},s_{k}],\
k=p+1,p+2,\dots, $ from  Lemma 4 for $\Delta=\R^n$  we obtain
\be{334}\begin{split}
&V(t,x^*(t))\\
&\leq V(t^*,x^*(t^*))-\frac{1}{\Gamma(q)}\Big(\int_{t^*}^{t_{p}} (t-s)^{q-1}c(||x^*(s)||)ds+\sum_{i=p}^{k-1}\int_{s_i}^{t_{i+1}} (t_{i+1}-s)^{q-1}c(||x^*(s)||)ds\Big)\\
&\leq V(t^*,x^*(t^*)). \nonumber
\end{split}
\ee
Therefore,  for $t\geq t^*$ the following inequality is satisfied:
\be{678} V(t, x^*(t))\leq V(t^*, x^*(t^*)).
\ee
 Then for
any $t\geq t^*$ applying (\ref{678}), condition 3(iii)  and inequality (\ref{88}) we get the inequalities \be{4427}
\begin{split}
b(||x^*(t)||)\leq V(t, x^*(t))\leq V(t^*, x^*(t^*))\leq a(||x^*(t^*)||)\leq a(\gg)<b(\epsilon).\nonumber
\end{split}
\ee   Therefore, the inequality  (\ref{4450}) holds for all $t\geq t^*$ (hence
for $t\geq t_0+T$).

\hfill~~~\sq

\begin{Rem} If the initial time $t_0$ is in an interval of noninstantaneous  impulses, i.e. $t_0\in \cup_{k=1}^{\infty}(t_k,s_k]$ then the results of Theorem 1, Theorem 2  and  Theorem 3  will be similar with slight changes in Definition 2 and  condition 2(ii) (Theorems 1,2) or condition 3(ii)(Theorem 3). 
\end{Rem}

\subsection{Stability by generalized Dini fractional derivative}\label{s6.2}

The proofs below are similar to those in Theorem 1 and Theorem 2
where Lemma 6 is used instead of Lemma 2.

 \begin{Th} ({\it Stability}).  Let the conditions of Theorem 1 be
 satisfied where the inequality (\ref{13517}) is replaced by
  \be{1352}\ \text{ }_{(\ref{1})}^c{\mathcal D}_+^qV(t,x)\leq0  \ \ \ \ \ \ \ \mbox{for}\ \  t\in (s_k,t_{k+1}],\ k=0,1,2,\dots,\ \  x\in \R^n.\ee
     Then  the zero solution of the NIFrDE (\ref{1}) is   stable.
     \end{Th}

  \begin{Th} ({\it Uniform stability}). Let the conditions
  of Theorem 2 be satisfied where the inequality (\ref{135}) is replaced by
\be{1353}\ \text{ }_{(\ref{1})}^c{\mathcal D}_+^qV(t,x)\leq 0  \ \ \ \ \ \ \ \mbox{for}\ \  t\in (s_k,t_{k+1}],\ k=0,1,2,\dots,\ \  x\in S(\gl).\ee
     Then  the zero solution of the IFrDE (\ref{1}) is  uniformly stable.
     \end{Th}

We will give the comparison result with the quadratic Lyapunov
function $V(x)=x^Tx=\sum_{k=1}^n x_k^2\in \Lambda^C(\R^n),\ \
x=(x_1,x_2,\dots,x_n)$.

\begin {Cor} ({\it Stability by a quadratic function}). Let condition (H) be satisfied and
    \begin{itemize}\item[(i)]  $x^T f(t,x)\leq 0$  for $t\in (s_k, t_{k+1}),\ \ k=0,1,2,\dots, \ \  x\in S(\gl)$;
\item[(ii)] for any $k=1,2,3,\dots$ and $x\in  S(\gl)$ for $\ t\in [s_k,t_k]$ the inequality $\big(\phi_k(t,x)\big)^T\phi_k(t,x)\leq x^Tx\ \ \mbox{for}\  t\in (t_k,s_k]$ holds;
\item[(iii)]  $b(||x||)\leq x^Tx\leq a(||x||)$ for $ t\in\R_+, \ x\in \R^n,$
  where  $a,b \in \mathcal{K}$.
\end{itemize}
     Then  the  zero solution of (\ref{1})  is  uniformly stable.
   \end {Cor}
	
	\begin {Rem} If $t_k=s_k$ and $y=\phi_k(t_k,y),\ y\in\R^n$ for all  $k=1,2\dots$  system (\ref{1})  reduces to a system of fractional differential equations. For the reduced (\ref{1}) in \cite {Yy}  the author defines a generalized
fractional-order derivative ( Dini-like derivative) in the Caputo sense  based on the fractional-order Dini derivative in the Caputo sense  (\cite {Lak fde}, \cite{LL}) and some sufficient conditions for  the stability with initial time difference are obtained.  
	
	\end{Rem}
	
\section{APPLICATIONS}\label{s7}

\subsection{Quadratic Lyapunov  function}\label{s7.1} We will apply the quadratic Lyaunov function and its fractional derivative.

\textbf{Example 6}.  Consider the scalar NIFrDE  (\ref{311}) where
$A\leq  0$ and $\phi_k(t,x)=a_k(t)x, \ a_k:\ [t_k,s_k]\to\R,  \
k=1,2,3,\dots$  are such that $\sup_{t\in[t_k,s_k]}|a_k(t)|\leq 1$.
Consider the quadratic Lyapunov function  $V(x)=x^2$. Then
$xf(t,x)=Ax^2\leq 0$  for $t\in (s_k, t_{k+1}),\ \ k=0,1,2,\dots, \
\  x\in S(\gl)$, i.e. condition (i) of Corollary 3 is satisfied.
Also, $\Big(\phi_k(t,x)\Big)^2=\Big(a_k(t)\Big)^2x^2 \leq x^2\ \
\mbox{for}\  t\in [t_k,s_k]$, i.e. condition (ii) is satisfied.

From Corollary 3 the zero solution of the scalar NIFrDE (\ref {311})
is uniformly stable (this was proved  directly in Example 2).

\hfill~~~\sq

\textbf{Example 7}.  Let the points $t_k, s_k, k=0,1,2,\dots, $ be
such that $0\leq t_k<s_k<t_{k+1},\ \lim_{k\to\infty}t_k=\infty$.
Consider the scalar noninstantaneous impulsive Caputo fractional
differential equation \be{347}\begin{split}
&\ _{s_k}^{c}D^{q}x=-a(t)x(1+x^2)\ \ \mbox{for} \ \ t\in(s_k,t_{k+1}], \ k=0,1,2,\dots,\\
&x(t)=c_k(t)x(t_k-0)\ \ \mbox{for} \ \ t\in(t_k,s_k], \ \ \ k=1,2,3,\dots, \\
&x(0)=x_0
\end{split}
\ee
 where $x \in \R$,   $a(t)\in C(\cup_{k=0}^\infty(s_k,t_{k+1}],\R_+)$,  $c_k(t)\in C([t_k,s_k],[-1,1])$, $k=1,2,\dots$.

 Consider the function $V(t,x)=x^2$.  Then $xf(t,x)=-a(t)x^2(1+x^2)\leq 0$, i.e. condition (i) of Corollary 3 is satisfied.
Also, $\Big(c_k(t)x\Big)^2=(c_k(t)x)^2\leq x^2\ \ \mbox{for}\  t\in
(s_k,t_k], \ \ \ k=1,2,3,\dots$, i.e. condition (ii) is satisfied.

 From Corollary 3, the trivial solution of NIFrDE (\ref{347}) is uniformly  stable.

 \hfill~~~\sq

\subsection{General Lyapunov  function}\label{s7.2} Now we will apply a general Lyapunov function and its generalized Caputo fractional Dini derivative.

\textbf{Example 8}. Let points $s_k= (4k+1)\frac{\pi}{2}, t_k=(4k-1)\frac{\pi}{2}, k=1,2,\dots$, $s_0=0$. Consider the following initial value problem for the scalar noninstantaneous impulsive
Caputo fractional differential equation {\footnotesize \noindent }%
\be{87}\begin{split}
&\ _{s_k}^{c}D^{0.5}x(t)=xf(t),\ \ \  t\in (s_k,t_{k+1}],\ k=0,1,2,\dots,\\
&x(t)=c_k(t) x(t_k-0),\ \ \ \ \ t\in [t_k,s_k],\ k=1,2,\dots,\\
& x({0})=x_{0},\end{split}\ee
where $x\in \R$, $c_k\in C([t_k,s_k],[-1,1])$, $f(t)=0.5\frac{\frac{-2}{\sqrt{t\pi}} +\sqrt{t}E_{2,1.5}(-t^2)}{2-sin(t)}$ , $k=0, 1, 2,\dots$.

Let $V(t,x)=x^2$.   Then $x\big(xf(t)\big)=x^2f(t)$.
Since the
sign of the function
$f(t)$ changes (see Figure 2) Corollary 3  and the  quadratic Lyapunov function  are not
applicable to the fractional equation (\ref{87}). 

Let $V(t,x)=(2-\sin(t))x^2$.

Apply generalized  Dini  fractional derivative given by (\ref{14}),  Example 4 and formula (\ref{16}) and we get
$$\text{
}_{(\ref{87})}^{c}{\mathcal D}_{+}^{0.5}V(t,x) =2 x^2\ (2-\sin(t))f(t).$$
The sign of the function $f(t)$ and the derivative $\text{
}_{(\ref{87})}^{c}{\mathcal D}_{+}^{0.5}V(t,x)$  are changeable. Therefore,  the application of fractional Dini derivative (\ref{14}) does not give us a conclusion about stability properties of NIFrDE (\ref {87}).

Now apply generalized Caputo fractional Dini derivative to the considered Lyapunov function. According to Example 4 for $t\in (s_k,t_{k+1}),\ k=0,1,2,\dots$ and $\Gamma({0.5})=\sqrt{\pi}$ we obtain
\be{1672}\begin{split}
&\text{
}_{(\ref{87})}^{c}D_{+}^{0.5}V(t,x;0, x_0 ) =2x^2(2-\sin(t))f(t) +x^2\
_{0}^{RL}D^{0.5}(2-\sin(t))
-\frac{2(x_{0})^2}{\sqrt{t\pi}}\\
& \ \ \ =2x^2(2-\sin(t))f(t) +x^2\Big(\frac{2}{\sqrt{t\pi}} -\sqrt{t}E_{2,1.5}(-t^2)\Big)
-\frac{2(x_{0})^2}{\sqrt{t \pi}}\leq 0.\end{split}\ee

Also, for $t\in [t_k,s_k],\ k=1,2,\dots$ we get
$V(t,c_k(t)x)=\Big(2-\sin(t)\Big)(c_k(t)x)^2\leq \Big(2-\sin(t)\Big)x^2\leq \Big(2-\sin(t_k)\Big)x^2=\Big(2-\sin((4k-1)\frac{\pi}{2})\Big)x^2= V(t_k-0,x)$, i.e. condition 2(ii) of Theorem 1 is satisfied.

According to Theorem 1 the zero solution of (\ref {87}) is  stable.

\hfill~~~\sq

\section{CONCLUSSIONS}\label{s8} 

Piecewise continuous scalar Lyapunov functions are applied to study
stability, uniform stability and asymptotic uniform stability of the
zero solution of nonlinear Caputo fractional differential equations
with not instantaneous impulses. Two types of derivatives of
Lyapunov function among NIFrDE are introduced and their applications
are discussed. Several sufficient conditions for stability, uniform
stability and uniform asymptotic  stability of the zero solution of
nonlinear NIFrDE are obtained. The results are illustrated with
several examples.

Note the above description and all consideration in  the paper could
be generalized when additionally an instantaneous jump at the points
$s_k$ are given and the right side part of the fractional equation is changing (so called variable structure),  i.e. problem (\ref{1}) can be generalized to
\be{199}\begin{split}
& \text{\ }_{s_k}^{c}D^{q}x=f_k(t,x) \ \mbox {for} \ t\in (s_k,t_{k+1}],\ k=0,1,2,\dots\\
& x(t)=\phi_i(t,  x(t_i-0))\ \ \ \mbox {for} \  t\in(t_i,s_i], \ i=1,2,\dots,\\
&x(s_i+0)=G_k(x(s_i-0) )\ \ \ \mbox {for} \  i=1,2,\dots,\\
& x(t_0)=x_0. \end{split}
\nonumber
\ee

\noindent{\sc Acknowledgments.} Research was partially supported by  the Fund NPD, Plovdiv University, No. MU15-FMIIT-008.


\begin{thebibliography}{39}

\bibitem {ABT}
S. Abbas, M. Benchohra, J. Trujillo, Upper and lower solution method for partial fractional differential inclusions  with not instantaneous impulses, {\it Progr. Fract. Differ. Appl.}, {\bf 11}, 1, (2015), 11--22.

\bibitem{AB}
R. Agarwal, M. Benchohra,  B. A. Slimani, Existence
results for differential equations with fractional order
and impulses, {\it Mem. Differ. Equ. Math. Phys.,} {\bf  44}, (2008), 1--21.

\bibitem{AH}
R. Agarwal, S. Hristova, Strict stability in terms of two measures for impulsive differential equations with ‘supremum’, {\it Appl. Anal. }, {\bf 91}, 7, (2012), 1379--1392.

\bibitem{AHR}
R. Agarwal, D. O'Regan, S. Hristova,  Stability of Caputo fractional differential equations by Lyapunov functions, {\it Appl. Math.} (accepted).

\bibitem{DM1}
 N. Aguila-Camacho,  M. A. Duarte-Mermoud, J. A. Gallegos,
 Lyapunov functions for
fractional order systems, {\it Comm. Nonlinear Sci. Numer.
Simul.}, {\bf 19}, (2014), 2951--2957.

\bibitem{H4}
D. Bainov, S. Hristova, The method of quasilinearization for the periodic boundary value problem for systems of impulsive differential equations, {\it Appl. Math. Comput.}, {\bf 117}, 1, (2001),  73--85. 

\bibitem{BK}
K. Balachandran, S. Kiruthika, Existence of solutions of abstract fractional impulsive semilinear evolution equations, {\it Electron. J. Qual. Theory Differ.
Equ.}, {\bf 4 }, (2010), 1--12.

\bibitem{BM}
D. Baleanu, O.G. Mustafa, On the global existence of solutions to a class of fractional differential
equations, {\it Comput. Math. Appl.,} {\bf59}, (2010), 1835--1841.

\bibitem{BS}
 M. Benchohra, D. Seba, Impulsive fractional differential equations in Banach spaces, {\it Electron. J. Qual. Theory Differ. Equ., } (Special
Edition I), {\bf 8},  (2009), 1--14.

\bibitem{B}
M. Benchohra, B. A. Slimani, Existence and uniqueness
of solutions to impulsive fractional differential equations,{\it
Electronic J. Differential Equations}, {\bf 2009},   10, (2009),   1--11.



\bibitem{1}
Sh. Das, {\it Functional Fractional Calculus}, Springer-Verlag Berlin Heidelberg, 2011.




\bibitem{devi vlm}  J. V. Devi , F.A. Mc Rae, Z. Drici, Variational
Lyapunov method for fractional differential equations, {\it Comput.
Math. Appl.}, {\bf 64},  (2012),  2982--2989.

\bibitem{2}
K. Diethelm,  {\it The Analysis of Fractional Differential
Equations}, Springer-Verlag Berlin Heidelberg,  2010.


\bibitem {FW}
M. Feckan, J.R. Wang, Y. Zhou, Periodic solutions for nonlinear evolution equations with non-istantaneous impulses, {\it Nonauton. Dyn. Syst.},  {\bf  1}, (2014), 93--101.

\bibitem{GD}
G. R. Gautam, J. Dabas, Mild solution for fractional functional integro-differential equation with
not instantaneous impulse, {\it Malaya J. Mat.}, {\bf 2},  3, (2014), 428–-437.


\bibitem {RD}
G. R. Guantam, J. Dabas, Mild solutions for a class of neutral fractional functional differential equations  with not instantaneous impulses, {\it Apl. Math. Comput.}, {\bf 259}, (2015), 480--489.

\bibitem {RD1}
G. R. Guantam, J. Dabas, Existence results on fractional functional differential equations  with not instantaneous impulses, {\it Int. J. Adv. Appl. Math. Mech.}, {\bf 1}, 3, (2014), 11--21.

\bibitem{HR}
E. Hernandez, D. O'Regan, On a new class of abstract impulsive differential equations, {\it  Proc. Amer. Math. Soc.}, {\bf 141},
(2013), 1641--1649.


\bibitem {S}
 S. Hristova, {\it Qualitative investigations and approximate methods for impulsive equations}, Nova Sci. Publ. Inc., New York, 2009.

\bibitem {H1}
 S. Hristova, Integral stability in terms of two measures for impulsive functional differential equations, {\it Math. Comput, Modell,}, {\bf 51},  1–2, (2010),  100–-108.



\bibitem {H2}
S. Hristova, Stability on a cone in terms of two measures for impulsive differential equations with “supremum”, {\it Appl. Math. Lett.}, {\bf 23}, 5, (2010), 508--511.

\bibitem{H3}
S. Hristova, Razumikhin method and cone valued Lyapunov functions for impulsive differential equations with “supremum”, {\it Inter. J. Dynam. Syst. Diff. Eq.}, {\bf 2}, 3-4, (2009), 223--236. 

\bibitem{H5}
S. Hristova, Lipschitz stability for impulsive differential equations with ‘supremum’, {\it International Electronic Journal of Pure and Applied Mathematics}, {\bf 1}, 4, (2010), 345--358.

\bibitem {BH1}
S. Hristova, D. Bainov, Periodic solutions of quasilinear nonautonomous systems with impulses, {\it Bull. Austral. Math. Soc.}, {\bf 31}, (1985), 185--198.

\bibitem {BH2}
S. Hristova, D. Bainov, Application of Liapunov's functions for studying the boundedness of solutions of systems with impulses, {\it COMPEL}, {\bf 5}, 1, (1986), 23--40.

\bibitem{HS}
S. Hristova, K. Stefanova, Practical stability of impulsive differential equations with “supremum” by integral inequalities, {\it Eur. J. Pure  Appl. Math.}, {\bf 5}, 1, (2012), 30--44. 



\bibitem{H}
J.B. Hu, G.P. Lu, S.B. Zhang, L.-D. Zhao, Lyapunov stability theorem
about fractional system without and with delay, {\it Commun.
Nonlinear Sci. Numer. Simul.}, {\bf 20}, (2015), 905--913.

\bibitem {E}
P. Eloe, S. Hristova, Method of the quasilinearization for nonlinear impulsive differential equations with linear boundary conditions, {\it EJQTDE}, {\bf 10}, (2002), 1-14. 

\bibitem{KP}
P. Kumar, D. N. Pandey, D. Bahuguna, On a new class of abstract impulsive functional
differential equations of fractional order, {\it J. Nonlinear Sci. Appl.}, {\bf 7},  (2014), 102--114.



\bibitem{LB}
V. Lakshmikantham, D.D. Bainov, P.S. Simeonov, {\it Theory of Impulsive Differential Equations}, World Scientific, Singapore, 1989.

\bibitem{Lak fde} V. Lakshmikantham,  S. Leela,  J.V. Devi, {\it Theory of
Fractional Dynamical Systems}, Cambridge Scientific Publishers, 2009.

\bibitem{LL}
V. Lakshmikantham,  S. Leela, M. Sambandham, Lyapunov theory for
fractional differential equations, {\it Commun. Appl. Anal.}, {\bf
12}, 4, (2008), 365--376.


\bibitem{Li1}
Y. Li, Y. Chen, I. Podlubny, Stability of fractional-order nonlinear
dynamic systems: Lyapunov direct method and generalized
Mittag-Leffler stability, {\it Comput. Math. Appl.}, {\bf 59}, (2010),
1810--1821.


\bibitem{LX}
P. Li, Ch. Xu, Boundary value problems of fractional order
differential equation with integral boundary conditions and not
instantaneous impulses, {\it J. Funct. Spaces}, {\bf 2015}, 2015,
Article ID 954925.



\bibitem {Li}
C. Li, D. Qian, Y. Chen, On Riemann-Liouville and Caputo Derivatives, {\it Discr. Dynam.
Nat, Soc.},
{\bf 2011}, Article ID 562494.


\bibitem{Li3}
C.P. Li,  F.R. Zhang, A survey on the stability of fractional
differential equations, {\it Eur. Phys. J., Special Topics}, {\bf
193},  (2011), 27--47.


\bibitem{LW1}
Y. M.  Liao, J. R.  Wang, A note on stability of impulsive differential equations, {\it Boundary Value Probl.}, {\bf 2014}, 2014:67.

\bibitem{LW2}
Z.Lin, J.R. Wang, W. Wei, Multipoint BVPs for generalized impulsive fractional differential equations, {\it Appl. Math. Comput.}, {\bf 258}, (2015), 608-616.

\bibitem{LW3}
Z.  Lin,  W. Wei, J. R. Wang, Existence and stability results for impulsive integro-differential equatons, {\it Facta Univer. (Nis)}, ser. Math. Inform., {\bf 29}, No 2 (2014), 119–130.



\bibitem{P}
D.N. Pandey, S. Das, N. Sukavanam, Existence of solutions for a second order neutral differential equation with state dependent delay and not instantaneous impulses, {\it Intern. J. Nonlinear Sci.}, {\bf 18}, 2, (2014), 145--155.

\bibitem {PRR}
M. Pierri, D. O'Regan, V. Rolnik, Existence of solutions for
semi-linear abstract differential equations with not instantaneous
impulses, {\it Appl. Math. Comput.}, {\bf 219}, (2013), 6743--6749.

\bibitem{podlubny} I. Podlubny, {\it Fractional Differential Equations}, Academic
Press, San Diego, 1999.



\bibitem{kilbas}  G. Samko,  A.A. Kilbas, O.I. Marichev, {\it Fractional Integrals
and Derivatives: Theory and Applications}, Gordon and Breach,  1993.


\bibitem{SS}
A. Sood,  S. K. Srivastava,  On Stability of Differential Systems with
Noninstantaneous Impulses, {\it Math. Probl. Eng.}, {\bf 2015}, 2015, Article ID 691687.


\bibitem{T}
J.C. Trigeassou, N. Maamri, J. Sabatier, A. Oustaloup, A Lyapunov
approach to the stability of fractional differential equations, {\it
Signal Processing}, {\bf 91}, (2011), 437--445.

\bibitem {ZN}
G. Wang, B. Ahmad, L. Zhang,  J. Nieto, Comments on the concept of existence of solution for impulsive
fractional differential equations, {\it Commun. Nonlinear Sci. Numer. Simulat.}, {\bf 19}, (2014),  401--403.





\bibitem {LW}
J. R. Wang, X. Li, Periodic BVP for integer/fractional order nonlinear differential equations with non-instantaneous impulses, {\it J. Appl. Math. Comput.}, {\bf 46},  1-2, (2014),  321--334.



\bibitem  {WL}
J. Wang, Z. Lin, A class of impulsive nonautonomous
differential equations and Ulam $-$ Hyers$-$Rassias
stability,
{\it Math. Meth. Appl. Sci}, {\bf 38},  5, (2015),  868--880.


\bibitem  {WZ}
Z. Wang, Y. Zang
Existence and stability of solutions to nonlinear impulsive differential equations in $\gb$-normed space, {\it Elect. J. Diff. Eq. } {\bf 2014},  83, (2014), 1--10.

\bibitem{Yy}
C. Yakar, Fractional differential equations in terms of comparison results and Lyapunov stability with initial time difference, {\it Abstr. Appl.  Anal.}, {\bf 2010}, 2010, Article ID 762857, Hindawi Publishing Corporation.


\bibitem{11}
Z. Yan, F. Lu, Existence results for a new class of fractional impulsive partial neutral stochastic integro-differential equations with ifinite delays, {\it J. Appl. Anal. Comput.}, {\bf 5}, 4, (2015), 329--346.


\bibitem{Y}
X. Yu, Existence and $\beta$-Ulam-Hyers stability for a class of fractional differential equations with non instantaneous  impulses, {\it Adv. Diff. Eq., } {\bf 2015}, (2015):104.\
\end{thebibliography}
\end{document}